\documentclass{article}

\usepackage[english]{babel}

\usepackage[letterpaper,top=2cm,bottom=2cm,left=3cm,right=3cm,marginparwidth=1.75cm]{geometry}
\usepackage{subcaption}
\usepackage{wrapfig}
\usepackage{amssymb}
\usepackage{amsmath}
\usepackage{graphicx}
\usepackage[colorlinks=true, allcolors=blue]{hyperref}
\usepackage[export]{adjustbox}
\title{Your Paper}
\author{You}
\begin{document}

\begin{center}
{\bf \large Inverse Source Problems for the Time-Fractional Evolution Equation}
\end{center}
\begin{center}
\textbf{A.A. Rahmonov\footnote{ V.I.Romanovskiy Institute of Mathematics,  Uzbekistan Academy of Sciences; \\ Bukhara State University, Bukhara, Uzbekistan; bsumath07@gmail.com. } 
}

\end{center}

{\bf Abstract.} In this paper, we investigate the direct and linear inverse problems of identifying time-dependent and time-independent source terms in a time-fractional diffusion-wave equation, using measured data at an interior point of the time interval. We first establish the existence, uniqueness, and regularity of the solution to the direct problem by employing the spectral method. Then, based on the properties of the direct problem, we study two inverse problems: one involving the recovery of a time-independent source term, and the other concerning the identification of a time-dependent coefficient function in the source term.
\vspace{0.3cm}

\vspace{0.3cm}


\section{Introduction}
Let $T$ be a positive constant and let $H$ be a separable Hilbert space with the scalar product $(\cdot,\cdot)$
and the norm $\|\cdot\|$,
 and $A:H\rightarrow H$
 be an arbitrary unbounded positive selfadjoint operator in $H$ and $A^{-1}$ is a compact operator. Suppose that $A$ has a complete in $H$ system of orthonormal eigenfunctions $e_n$
 and a countable set of positive eigenvalues $\lambda_n$. It is convenient to assume that the eigenvalues do not decrease as their number increases, i.e., $0<\lambda_1\le\lambda_2\le \dots\rightarrow+\infty$. We consider the following fractional diffusion-wave equation:
\begin{equation}\label{Eq1}
\partial_{t}^{\alpha}u(t)+Au(t)=F(t),\quad t\in(0,T),    
\end{equation}
with the Gerasimov-Caputo time fractional derivative $\partial_{t}^{\alpha}$ of order $1<\alpha<2$, defined by
\begin{equation*}
\partial_{t}^{\alpha}y(t)=\frac{t_+^{1-\alpha}}{\Gamma(2-\alpha)}\star y''(t),
\end{equation*}
where $\star$ denotes the convolution, $\Gamma(\cdot)$ is Euler's gamma function and
\begin{equation*}
    t_{+}^{1-\alpha}=\begin{cases}
        t^{1-\alpha},&t>0,\\
        0,&t\le0,
    \end{cases}
\end{equation*}  
and the source term $F(t)$ is given.

We supplement the above fractional wave equation with the non-local initial conditions
\begin{equation}\label{Eq2}
    \gamma u(0)+u(T)=\varphi,\quad \partial_tu(T)=\psi,
\end{equation}
here, we assume that the $\gamma$ is given number, $\varphi$ and $\psi$ are given elements of $H$.

The fractional diffusion-wave equation with Gerasimov-Caputo time derivatives of order $1<\alpha<2$ has attracted substantial interest due to its ability to model phenomena exhibiting both diffusive and wave-like characteristics.  Foundational work by Sakamoto and Yamamoto \cite{Sakamoto2011} rigorously established the well-posedness and regularity of initial and boundary value problems for such equations, alongside important results on related inverse problems including source and coefficient identification. Building on this foundation, recent studies have extended the analysis to fractional PDEs with nonlocal or multi-point initial conditions and have explored backward and ill-posed inverse problems. For example, Alimov and Ashurov \cite{Alimov2020} examined backward problems for time-fractional subdiffusion equations, developing stability estimates and regularization methods critical to handling the inherent ill-posedness of inverse problems in fractional contexts. Similarly,  Ashurov and Fayziev \cite{Ashurov2021} studied the fractional diffusion equation with nonlocal initial conditions, providing existence, uniqueness, and stability results for direct and inverse source problems within spectral frameworks. Moreover, Floridia, Li and Yamamoto \cite{Florida2020} analyzed the well-posedness of backward problems in time for general time-fractional diffusion equations, contributing important theoretical insights into stability and uniqueness for such inverse problems. Such nonlocal-in-time conditions arise naturally in applications where control or observation is performed at the terminal time or across time intervals, and have been rigorously studied in recent literature. In particular, Ruzhansky et al. \cite{Ruzhansky2018} developed a detailed spectral framework for fractional diffusion-wave equations under multi-point and integral conditions, proving existence and uniqueness via Mittag-Leffler expansions. Their approach is particularly well-suited to Hilbert space settings, where the operator 
$A$ admits a complete set of eigenfunctions.  These advancements motivate the present study, where we consider fractional diffusion-wave equations in a separable Hilbert space equipped with a positive self-adjoint operator and impose nonlocal time conditions. Our work extends the existing theory by addressing novel inverse problems involving time-independent sources and time-dependent coefficients, utilizing spectral methods and fractional calculus to establish uniqueness, stability, and constructive solution representations.

We begin by analyzing the direct problem, establishing existence, uniqueness, and regularity of the solution in appropriate function spaces. The solution is represented using eigenfunction expansions and the properties of the two-parameter Mittag-Leffler function. This representation forms the basis for the study of two inverse problems. 

In Inverse Problem 1, we assume that the right-hand side is time-independent, i.e., $F(t)\equiv f\in H$ and seek to determine both the function $u(t)$ and the unknown source term $f,$ given additional overdetermination data at an intermediate time, i.e.,
\begin{equation}\label{Eq3}
  u(\xi)=h,\quad 0<\xi<T,  
\end{equation}
where $h\in D(A)\subset H$ is a given element. This type of problem is motivated by applications in source identification and has been studied in various forms in \cite{Kirane2012}-\cite{Lingyun2024}. These works employ both spectral and analytic continuation techniques to establish uniqueness and stability under suitable assumptions.

{\bf Definition 1.1.} A pair $\{u(t),f\}$ of functions $u\in C([0,T];D(A))\cap C^1([0,T];H)$ and $f\in H$ with the properties $\partial_t^{\alpha}u(t),Au(t)\in C((0,T];H)$ and satisfying conditions \eqref{Eq1}-\eqref{Eq3} is called the solution of the {\it inverse problem 1}.

In Inverse Problem 2, we consider the identification of a time-dependent coefficient $p(t)$ in the source term $F(t)=p(t)f,$ where $f\in H$ is known, and the overdetermination condition involves a scalar measurement   
\begin{equation}\label{Eq4}
  \Phi[u(t)]=h(t),\quad 0<t<T, 
\end{equation}
where $h:[0,T]\rightarrow\mathbb{R}$ is a given function, $\Phi:D(\Phi)\subset H\rightarrow\mathbb{R}$ is a known linear bounded functional, where $D(\Phi)=\{u\in H: u\in H\}$. This formulation arises in applications where internal parameters change over time and are tracked through observable functionals. Recent studies by \cite{Durdiev2023} and \cite{Turdiev2023} have addressed similar coefficient identification problems using Volterra integral equations and compactness arguments. While related inverse source and coefficient problems have been studied in fractional diffusion settings  \cite{Kirane2012} and \cite{Xing2023}, the present work extends the analysis to fractional diffusion-wave equations with nonlocal time conditions and general self-adjoint operators. Our results include new existence, uniqueness, and reconstruction strategies in both inverse problems, with explicit spectral representations that differ significantly from classical approaches.

{\bf Definition 1.2.} A pair $\{u(t),p(t)\}$ of functions $u\in  C([0,T];D(A))\cap C^1([0,T];H)$ and $p\in L^{\infty}(0,T)$ with the properties $\partial_t^{\alpha}u(t),Au(t)\in C((0,T];H)$ and satisfying conditions \eqref{Eq1}, \eqref{Eq2}, \eqref{Eq4} is called the solution of the {\it inverse problem 2}.

This paper is organized as follows. In Section 2, we study the direct problem in detail, including the existence, uniqueness, and regularity of the solution. Section 3 is devoted to Inverse Problem 1, where we prove the uniqueness of the unknown time-independent source term and derive a constructive method for its reconstruction.  In Section 4, we consider Inverse Problem 2, focusing on the identification of a time-dependent coefficient function.

\section{Notations and preliminaries}

\subsection{Functional spaces.} 
\ \ \ \ In this subsection, we start by giving some definitions and preliminary results which will be used throughout this paper. 

Let $\sigma$ be an arbitrary real number. We introduce the power of operator $A$, acting in $H$ according to the rule  (see, e.g., \cite{Pazy} or \cite{Alimov1972})
\begin{equation*}
A^{\sigma}g=\sum\limits_{n=1}^{\infty}\lambda_n^{\sigma}g_ne_n,
\end{equation*}
where $g_n$ is the Fourier coefficients of a function $g\in H:\,g_n=(g,e_n)$. The domain of this operator has the form 
\begin{equation*}
    D(A^{\sigma})=\{g\in H:\,\sum\limits_{n=1}^{\infty}\lambda_n^{2\sigma}(g,e_n)^2<\infty\}.
\end{equation*}
For elements of $D(A^{\sigma})$ we introduce the norm
\begin{equation*}
\|g\|_{D(A^{\sigma})}^2=\sum\limits_{n=1}^{\infty}\lambda_n^{2\sigma}g_n^2=\|A^{\sigma}g\|^2.
\end{equation*}
Let $V$ be a Banach space, and we denote by $C([0,T];V)$ the space of all continuous functions from $[0,T]$ to $V$ endowed with the norm $\|g\|_{C([0,T];V)}=\max\limits_{t\in[0,T]}\|g(t)\|_V$, and by $C^{\theta}([0,T];V),\,0<\theta<1$, the subspace of $C([0,T];V)$ which includes all H\"{o}lder-continuous functions, and is equipped with the $\theta^{th}-$H\"{o}lder seminorm
\begin{equation*}
    [g]_{C^{\theta}([0,T];V)}:=\sup\limits_{0\le t_1<t_2\le T}\frac{\|g(t_2)-g(t_1)\|_V}{|t_2-t_1|^{\theta}},
\end{equation*}
and the $\theta^{th}-$H\"{o}lder norm
\begin{equation*}
    \|g\|_{C^{\theta}([0,T];V)}=\|g\|_{C([0,T];V)}+[g]_{C^{\theta}([0,T];V)}.
\end{equation*}
Moreover, we denote by $AC[0,T]$ the space of absolutely continuous functions on $[0,T]$, and define 
\begin{equation*}
    AC^m[0,T]:=\{y(t):\,y\in C^{m-1}[0,T],\,y^{(m)}(t)\in AC[0,T]\},\quad m\ge2. 
\end{equation*}

Now, we introduce the Mittag-Leffler function (see \cite{Kilbas}, pp.40-45)
$$
E_{\rho,\mu}(z)=\sum\limits_{k=0}^{\infty}\frac{z^k}{\Gamma( \rho k+\mu)},\quad z\in\mathbf{C}
$$
with $\mathrm{Re}(\rho)>0$ and $\mu\in\mathbf{C}$. It is known that  $E_{\rho,\mu}(z)$ is an entire function in $z\in \mathbf{C}$.

{\bf Lemma 1.1.} (\cite{Podlubny}, p.35) {\it Let $0<\rho<2$ and $\mu\in\mathbf{R}$ be arbitrary and $\theta$ satisfy $\frac{\pi\rho}{2}<\theta<\min\{\pi,\pi\rho\}$. Then there exists a constant $c=c(\rho,\mu,\theta)>0$ such that
$$
\left|E_{\rho,\mu}(z)\right|\leq\frac{c}{1+|z|},\qquad \theta\le |\mbox{arg} (z)|\le \pi.
$$
}

For the proof, we refer to \cite{Djrbashian}.

{\bf Lemma 1.2.} (\cite{Liao}) {\it Let $1<\alpha<2$ and $h(t)\in AC[0,T]$. Define
$$
p(t)=h\star t^{\alpha-1}E_{\alpha,\alpha}(-\lambda t^{\alpha}),\quad t>0.
$$
If $\lambda>0$, then $p(t)\in AC^2[0,T]$ satisfies
\begin{equation*}
    \partial_t^{\alpha}p(t)+\lambda p(t)=h(t),\quad 0<t\le T.
\end{equation*}
}

{\bf Lemma 1.3.} {\it Let $1<\alpha<2$ and $\lambda>0$, then we have
$$
\begin{cases}
    \frac{d}{dt}E_{\alpha,1}(-\lambda t^{\alpha})=-\lambda t^{\alpha-1}E_{\alpha,\alpha}(-\lambda t^{\alpha}),&t\ge0,\vspace{0.2cm}\\
      \frac{d}{dt}tE_{\alpha,2}(-\lambda t^{\alpha})=E_{\alpha,1}(-\lambda t^{\alpha}),&t\ge0,\vspace{0.2cm}\\
      \frac{d}{dt}t^{\alpha-1}E_{\alpha,\alpha}(-\lambda t^{\alpha})= t^{\alpha-2}E_{\alpha,\alpha-1}(-\lambda t^{\alpha}),&t>0,\vspace{0.2cm}\\
       \partial_t^{\alpha}(tE_{\alpha,2}(-\mu t^{\alpha}))=-\lambda tE_{\alpha,2}(-\lambda t^{\alpha}),&t>0,\vspace{0.2cm}\\
        \partial_t^{\alpha}E_{\alpha,1}(-\lambda t^{\alpha})=-\lambda E_{\alpha,1}(-\lambda t^{\alpha}),&t>0.
\end{cases}
$$
}
Please refer to \cite{Kilbas} for proof.

\subsection{Solution of the OFVP}\label{subsection:2.2}

Now, we discuss solutions of the FVP for the fractional ordinary differential equation
\begin{equation}\label{Eq5}
    \partial_t^{\alpha}y(t)+\lambda y(t)=f(t),\quad 0<t<T,
\end{equation}
and
\begin{equation}\label{Eq6}
    \gamma y(0)+y(T)=a_T,\quad y'(T)=b_T,
\end{equation}
where $\lambda,\,\gamma,a_T$ and $b_T$ are given real numbers. Moreover, we have assumed the $\lambda>0$ and $\gamma\not=0$.

According to Theorem 5.15 (see \cite{Kilbas}, p.323), the general solution to Eq. \eqref{Eq5} is given by
\begin{equation}\label{Eq7}
    y(t)=E_{\alpha,1}(-\lambda t^{\alpha})\cdot d_1+tE_{\alpha,2}(-\lambda t^{\alpha})\cdot d_2+f(t)\star t^{\alpha-1}E_{\alpha,\alpha}(-\lambda t^{\alpha}),
\end{equation}
where $d_i,i=1,2$ are arbitrary real constants. By Lemma 1.2 for $\alpha>0,\,t>0$, from \eqref{Eq7},  one has
\begin{equation}\label{Eq9}
    y'(t)=-\lambda t^{\alpha-1}E_{\alpha,\alpha}(-\lambda t^{\alpha})\cdot d_1+E_{\alpha,1}(-\lambda t^{\alpha}) \cdot d_2
+f(t)\star t^{\alpha-2}E_{\alpha,\alpha-1}(-\lambda t^{\alpha}).   
\end{equation}
Let $t=T$ in \eqref{Eq7} and \eqref{Eq9}, then we obtain
\begin{equation}\label{Eq10}
    y(T)=E_{\alpha,1}(-\lambda T^{\alpha})\cdot d_1+TE_{\alpha,2}(-\lambda T^{\alpha}) \cdot d_2+f(r)\star r^{\alpha-1}E_{\alpha,\alpha}(-\lambda r^{\alpha})\Big|_{r=T},
\end{equation}
and
\begin{equation}\label{Eq11}
    y'(T)=-\lambda T^{\alpha-1}E_{\alpha,\alpha}(-\lambda T^{\alpha})\cdot d_1+E_{\alpha,1}(-\lambda T^{\alpha}) \cdot d_2+f(r)\star r^{\alpha-2}E_{\alpha,\alpha-1}(-\lambda r^{\alpha})\Big|_{r=T}.
\end{equation}
Substituting \eqref{Eq10} and \eqref{Eq11} into the boundary conditions \eqref{Eq6}, one has
\begin{equation}\label{Eq12}
\begin{cases}
d_1=\frac{1}{\rho(\lambda T^{\alpha})}\Big[E_{\alpha,1}(-\lambda T^{\alpha})a_T-TE_{\alpha,2}(-\lambda T^{\alpha})b_T&\vspace{0.2cm}\\
\qquad\qquad\qquad-E_{\alpha,1}(-\lambda T^{\alpha})f(r)\star r^{\alpha-1}E_{\alpha,\alpha}(-\lambda r^{\alpha})\Big|_{r=T}&\vspace{0.2cm}\\
\qquad\qquad\qquad+TE_{\alpha,2}(-\lambda T^{\alpha})f(r)\star r^{\alpha-2}E_{\alpha,\alpha-1}(-\lambda r^{\alpha})\Big|_{r=T}\Big],\vspace{0.3cm}\\
d_2=\frac{1}{\rho(\lambda T^{\alpha})}\Big[(\gamma+E_{\alpha,1}(-\lambda T^{\alpha}))b_T+\lambda T^{\alpha-1}E_{\alpha,\alpha}(-\lambda T^{\alpha})a_T&\vspace{0.2cm}\\
\qquad\qquad\qquad-(\gamma+E_{\alpha,1}(-\lambda T^{\alpha}))f(r)\star r^{\alpha-2}E_{\alpha,\alpha-1}(-\lambda r^{\alpha})\Big|_{r=T}&\vspace{0.2cm}\\
\qquad\qquad\qquad\qquad-\lambda T^{\alpha-1}E_{\alpha,\alpha}(-\lambda T^{\alpha})f(r)\star r^{\alpha-1}E_{\alpha,\alpha}(-\lambda r^{\alpha})\Big|_{r=T}\Big],\vspace{0.3cm}
\end{cases}
\end{equation}
where 
\begin{equation}\label{Eq13}
    \rho(\eta)=\gamma E_{\alpha,1}(-\eta)+E_{\alpha,1}(-\eta)^2+\eta E_{\alpha,2}(-\eta)E_{\alpha,\alpha}(-\eta),\quad \eta>0.
\end{equation}
Putting the constants $d_1,d_2$ into \eqref{Eq7}, we arrive
\begin{equation*}
    y(t)=\frac{E_{\alpha,1}(-\lambda t^{\alpha})}{\rho(\lambda T^{\alpha})}\Big[E_{\alpha,1}(-\lambda T^{\alpha})a_T-TE_{\alpha,2}(-\lambda T^{\alpha})b_T\qquad\qquad\qquad\qquad\qquad\qquad\qquad\qquad\qquad\qquad
    \end{equation*}
    \begin{equation*}
-E_{\alpha,1}(-\lambda T^{\alpha})f(r)\star r^{\alpha-1}E_{\alpha,\alpha}(-\lambda r^{\alpha})\Big|_{r=T}\qquad\qquad\qquad\qquad\qquad\qquad\qquad\quad\qquad\,
  \end{equation*}
    \begin{equation*}
+TE_{\alpha,2}(-\lambda T^{\alpha})f(r)\star r^{\alpha-2}E_{\alpha,\alpha-1}(-\lambda r^{\alpha})\Big|_{r=T}\Big]\qquad\qquad\qquad\qquad\qquad\qquad\qquad\quad
\end{equation*}
\begin{equation*}
    +\frac{tE_{\alpha,2}(-\lambda t^{\alpha})}{\rho(\lambda T^{\alpha})}\Big[(\gamma+E_{\alpha,1}(-\lambda T^{\alpha}))b_T+\lambda T^{\alpha-1}E_{\alpha,\alpha}(-\lambda T^{\alpha})a_T\qquad\qquad\qquad\qquad\qquad\quad
\end{equation*}
\begin{equation*}
-(\gamma+E_{\alpha,1}(-\lambda T^{\alpha}))f(r)\star r^{\alpha-2}E_{\alpha,\alpha-1}(-\lambda r^{\alpha})\Big|_{r=T}\qquad\qquad\qquad\qquad\qquad\quad\quad\quad\,\,\,\,
\end{equation*}
\begin{equation}\label{Eq14}
\quad-\lambda T^{\alpha-1}E_{\alpha,\alpha}(-\lambda T^{\alpha})f(r)\star r^{\alpha-1}E_{\alpha,\alpha}(-\lambda r^{\alpha})\Big|_{r=T}\Big]+f(t)\star t^{\alpha-1}E_{\alpha,\alpha}(-\lambda t^{\alpha}).
\end{equation}

Before studying the problem \eqref{Eq1}-\eqref{Eq2}, we analyse the zeros of the function $\rho(\eta)$ defined by \eqref{Eq13}.

{\bf Lemma 2.1.} \textit{For $\gamma\not=0$, the set $\mathrm{P}=\{\eta>0:\,\rho(\eta)=0\}$ is a non-empty and finite set. Moreover, there exists a large constant $c_1>0$ such that 
\begin{equation}\label{Eq14.1}
    |\rho(\eta)|\ge\frac{c_1}{\eta},\quad\mbox{as}\quad \eta \ge c_1. 
\end{equation}
}

{\bf Proof.} \textit{$1^{st}$ case $\gamma<-1$.} Indeed, by the analyticity of the Mittag-Leffler function, we see that $\rho(\eta)$ is analytic in $\eta>0$ and continuous in $[0,\infty)$. Moreover, by the asymptotics of the Mittag-Leffler functions for $1<\alpha<2$, we see that
\begin{equation}\label{Eq15.0}
    \begin{cases}
        E_{\alpha,1}(-\eta)=\frac{1}{\Gamma(1-\alpha)}\frac{1}{\eta}+O\left(\frac{1}{\eta^2}\right),&         E_{\alpha,2}(-\eta)=\frac{1}{\Gamma(2-\alpha)}\frac{1}{\eta}+O\left(\frac{1}{\eta^2}\right),\\
                E_{\alpha,\alpha}(-\eta)=\frac{-1}{\Gamma(-\alpha)}\frac{1}{\eta^2}+O\left(\frac{1}{\eta^3}\right),&\mbox{as}\quad \eta\rightarrow\infty,
    \end{cases}
\end{equation}
(see, for example, \cite{Podlubny}, pp. 29-37). Therefore
\begin{equation*}
      \rho(\eta)= \gamma E_{\alpha,1}(-\eta)+E_{\alpha,1}(-\eta)^2+\eta E_{\alpha,2}(-\eta)E_{\alpha,\alpha}(-\eta)\qquad\qquad\qquad\qquad\qquad\qquad
\end{equation*}
\begin{equation*}
 =\frac{\gamma}{\Gamma(1-\alpha)}\frac{1}{\eta}+O\left(\frac{1}{\eta^2}\right)+\left(\frac{1}{\Gamma(1-\alpha)}\frac{1}{\eta}+O\left(\frac{1}{\eta^2}\right)\right)^2\qquad\qquad\qquad\quad\,\,
 \end{equation*}
\begin{equation*}
 -\eta\left(\frac{1}{\Gamma(2-\alpha)}\frac{1}{\eta}+O\left(\frac{1}{\eta^2}\right)\right)\left(\frac{1}{\Gamma(-\alpha)}\frac{1}{\eta^2}+O\left(\frac{1}{\eta^3}\right)\right)\qquad\qquad\qquad\quad\,\,
\end{equation*}
\begin{equation*}
 =\frac{\gamma}{\Gamma(1-\alpha)}\frac{1}{\eta}+O\left(\frac{1}{\eta^2}\right)\quad \mbox{as}\quad \eta\rightarrow\infty.\qquad\qquad\qquad\qquad\qquad\qquad\quad\,\,
\end{equation*}
Since $\Gamma(1-\alpha)<0$, and by our assumption $\gamma<-1$, we see that $\frac{\gamma}{\Gamma(1-\alpha)}>0$. Then, there exists a constant $M>0$ such that $\rho(\eta)>0$ for $\eta\ge M$. Since $\rho(0)=\gamma+1,$ by the continuity of $\rho$ in $[0,\infty)$, we can find a sufficiently small constant $\varepsilon>0$ such that $\rho(\eta)<0$ for $0\le\eta\le \varepsilon$. Therefore the intermediate value theorem yields that there exists $\eta_0\in (\varepsilon,M)$ such that $\rho(\eta_0)=0.$ Moreover, since $\rho$ is analytic in $[\varepsilon,M]$, the set $\{\eta\in[\varepsilon,M]:\,\rho(\eta)=0\}$ is a finite set. Otherwise $\rho(\eta)=0$ for each $\eta\in[\varepsilon,M]$, which implies $\rho(0)=0$ by the continuity of $\rho(\eta)$ at $\eta=0$, which contradicts $\rho(0)=\gamma+1$.

\textit{$2^{nd}$ case $\gamma>0$.} In this case, we see that $\frac{\gamma}{\Gamma(1-\alpha)}<0.$ Therefore, the same as above, the set $\mathrm{P}$ is $\not=\emptyset$ and finite. Because of $\rho(0)=\gamma+1>0.$

\textit{$3^{rd}$ case $\gamma=-1$.} Easy to see that $\rho(0)=0$ and $\frac{\gamma}{\Gamma(1-\alpha)}>0.$ Considering all this,  we don't have any information about the zeros of the function \eqref{Eq13}. If a function \eqref{Eq13} has zeros, then they are finite because it is analytic. To do this, we calculate the integral of the function \eqref{Eq13} in a finite interval $T\in[0, M]$ and study its behaviour for  $\lambda\gg1$.

By \eqref{Eq13} and \eqref{Eq15.0}, we obtain
\begin{equation*}
    \int_0^M\rho(\lambda T^{\alpha})dT=\frac{-M^{1-\alpha}}{\Gamma(2-\alpha)}\frac{1}{\lambda}+O\Big(\frac{1}{\lambda^2}\Big),\quad \mbox{as}\quad \lambda\gg1.
\end{equation*}

\begin{wrapfigure}{r}{0.55\textwidth}
\includegraphics[width=0.9\linewidth]{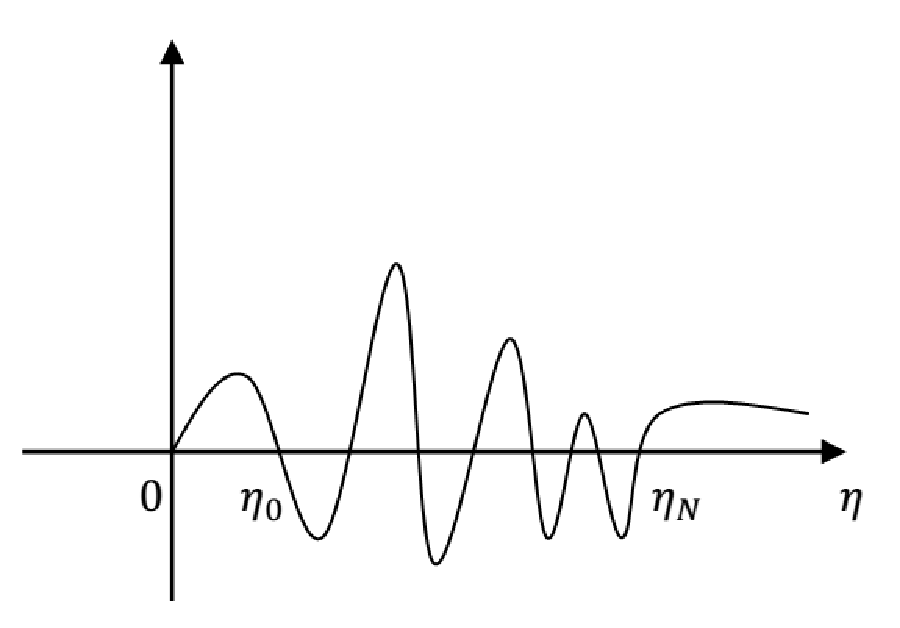} 
\caption{$\eta_0,..,\eta_N$ are the roots of $\rho(\eta)$}
\label{fig_1}
\end{wrapfigure}

The main term is $\frac{-M^{1-\alpha}}{\Gamma(2-\alpha)}<0$, so the function \eqref{Eq13} has zero (see Fig. \eqref{fig_1}).

\textit{$4^{th}$ case $\gamma\in(-1,0)$.} For the same reason as in case 3, the set $\mathrm{P}$ is not empty and finite.  

Furthermore, by the asymptotic expansions of Mittag-Leffler functions (see \eqref{Eq15.0}) with large $\eta$ there exists $c_1>0$ such that:
\begin{equation}\label{Eq15}
    |\rho(\eta)|\ge\frac{c_1}{\eta}. 
\end{equation}

Thus the proof of Lemma 2.1 is completed.
\\

We set 
\begin{equation}\label{Eq16}
    \{\eta_0,...,\eta_N\}=\{\eta>0:\,\rho(\eta)=0\}
\end{equation}
with $\eta_0<...<\eta_N$. By Lemma 2.1, we have no information on the number $N$ of the zeros of $\rho$, except that it exists. In the next lemma, we will provide an upper bound of the largest zero $\eta_N$.

aa

For convenience we set
\begin{equation*}
    \mu_1:=\frac{1}{\alpha\Gamma(-\alpha)},\quad \mu_2:=\frac{1}{\alpha(\alpha-1)\Gamma(-\alpha)},\quad \mu_3:=\frac{1}{\Gamma(-\alpha)}. 
\end{equation*}
Then, we have $\mu_1,\mu_2,\mu_3>0$ by $\Gamma(-\alpha)>0$ for $1<\alpha<2$. For $1<\alpha<2$, we choose $\theta$ such that 
\begin{equation}\label{Eq17}
    \frac{\pi\alpha}{2}<\theta<\pi.
\end{equation}

\begin{wrapfigure}{r}{0.5\textwidth}
\includegraphics[width=0.8\linewidth]{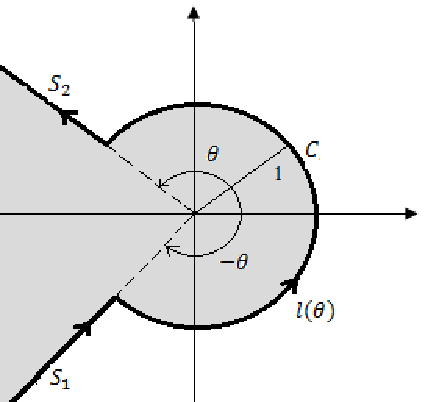} 
\caption{The contour of integration $l(\theta)$}
\label{fig_2}
\end{wrapfigure}

By $l=l(\theta)$ (see Fig.\eqref{fig_2}.) we denote the counter in $\mathbb{C}$ which is directed from $\infty e^{-\sqrt{-1}\theta}$ to $\infty e^{\sqrt{-1}\theta}$ and consists of (see \cite{Djrbashian}, p.126)
\begin{equation*}
    l(\theta)=\begin{cases}
        S_1=\{z:\,\mathrm{arg}\ z=-\theta,\,|z|\ge1\},\\
        C\,\,=\{z:\,-\theta\le\mathrm{arg}\ z\le\theta,\,|z|=1\},\\
        S_2=\{z:\,\mathrm{arg}\ z=\theta,\,|z|\ge1\}.
    \end{cases}
\end{equation*}
\vspace{2.7cm}

Moreover, we set
\begin{equation*}
    \nu_1=\frac{1}{2\pi\alpha \sin{\theta}}\int_{l}|\exp{(\zeta^{\frac{1}{\alpha}})}||\zeta||d\zeta|;\quad \nu_2=\frac{1}{2\pi\alpha \sin{\theta}}\int_{l}|\exp{(\zeta^{\frac{1}{\alpha}})}||\zeta^{1-\frac{1}{\alpha}}||d\zeta|;
\end{equation*}
\begin{equation*}
    \nu_3=\frac{1}{2\pi\alpha \sin{\theta}}\int_{l}|\exp{(\zeta^{1+\frac{1}{\alpha}})}||\zeta||d\zeta|.
\end{equation*}
Since there exists a constant $c_0>0$ such that $\exp{(\zeta^{\frac{1}{\alpha}})}\le \exp{(-c_0|\zeta|^{\frac{1}{\alpha}})}$ for any $\zeta\in l$ (see \cite{Djrbashian}, p.135), we can directly verify that $0<\nu_1,\nu_2,\nu_3<\infty$.

We are now able to prove the following lemma.

{\bf Lemma 2.2.} \textit{Let $\gamma>0$. Then, we have
\begin{equation*}
    \eta_N<\max\left\{\frac{1}{|\cos{\theta}|},\frac{1}{\mu_1\gamma}\left(\frac{1}{\alpha^2(\alpha-1)\Gamma(-\alpha)}+\gamma\nu_1+2\mu_1\nu_1+\mu_2\nu_3+\mu_3\nu_2+\nu_1^2+\nu_2\nu_3\right)\right\}.
\end{equation*}
}

{\bf Proof.} By formula (1.145) in \cite{Podlubny}, we have
\begin{equation}\label{Eq18}
    \begin{cases}
        E_{\alpha,1}(-\eta)=-\frac{\mu_1}{\eta}+I_{\alpha,1}(\eta),&         E_{\alpha,2}(-\eta)=\frac{\mu_2}{\eta}+I_{\alpha,2}(\eta),\\
                E_{\alpha,\alpha}(-\eta)=-\frac{\mu_3}{\eta^2}+I_{\alpha,\alpha}(\eta),&\eta\ge1,
    \end{cases}
\end{equation}
where
\begin{equation*}
    I_{\alpha,j}(\eta)=\frac{-1}{2\pi\alpha\sqrt{-1}\eta}\int_{l}\exp{(\zeta^{\frac{1}{\alpha}})}\zeta^{\frac{1-j}{\alpha}+1}\frac{d\zeta}{\zeta+\eta},\quad j=1,2,
\end{equation*}
\begin{equation*}
    I_{\alpha,\alpha}(\eta)=\frac{1}{2\pi\alpha\sqrt{-1}\eta^2}\int_{l}\exp{(\zeta^{\frac{1}{\alpha}})}\zeta^{\frac{1}{\alpha}+1}\frac{d\zeta}{\zeta+\eta},\quad \eta\ge1.
\end{equation*}
By \cite{Floridia}, we see that
\begin{equation}\label{Eq19}
    |I_{\alpha,1}(\eta)|\le\frac{\nu_1}{\eta^2},\quad |I_{\alpha,2}(\eta)|\le\frac{\nu_2}{\eta^2},\quad |I_{\alpha,\alpha}(\eta)|\le\frac{\nu_3}{\eta^3},\quad \mbox{for}\quad \eta\ge\frac{1}{|\cos{\theta}|}.
\end{equation}
Further, applying  \eqref{Eq18} and \eqref{Eq19} in \eqref{Eq13}, we obtain
\begin{equation*}
    \rho(\eta)=\gamma E_{\alpha,1}(-\eta)+E_{\alpha,1}(-\eta)^2+\eta E_{\alpha,2}(-\eta)E_{\alpha,\alpha}(-\eta)\qquad\qquad\qquad\qquad\qquad\qquad\quad\quad\,\,\,
\end{equation*}
\begin{equation*}
    =-\frac{\mu_1\gamma}{\eta}+\gamma I_{\alpha,1}(\eta)+\left(\frac{\mu_1}{\eta}-I_{\alpha,1}(\eta)\right)^2+\left(\frac{\mu_2}{\eta}+I_{\alpha,2}(\eta)\right)\left(-\frac{\mu_3}{\eta}+\eta I_{\alpha,\alpha}(\eta)\right)
\end{equation*}
\begin{equation*}
    =-\frac{\mu_1\gamma}{\eta}+\frac{1}{\eta^2}(\mu_1^2-\mu_2\mu_3)+\Big[(\gamma-2\frac{\mu_1}{\eta})I_{\alpha,1}(\eta)+I_{\alpha,1}(\eta)^2+\mu_2I_{\alpha,\alpha}(\eta)\qquad\quad\,\,\,
\end{equation*}
\begin{equation*}
    \,\quad-\frac{\mu_3I_{\alpha,2}(\eta)}{\eta}+\eta I_{\alpha,2}(\eta)I_{\alpha,\alpha}(\eta)\Big]\le-\frac{\mu_1\gamma}{\eta}+\frac{1}{\eta^2}|\mu_1^2-\mu_2\mu_3|+\Big[(\gamma+2\frac{\mu_1}{\eta})|I_{\alpha,1}(\eta)|
\end{equation*}
\begin{equation*}
    +I_{\alpha,1}(\eta)^2+\mu_2|I_{\alpha,\alpha}(\eta)|+\frac{\mu_3|I_{\alpha,2}(\eta)|}{\eta}+\eta |I_{\alpha,2}(\eta)||I_{\alpha,\alpha}(\eta)|\Big]\le-\frac{\mu_1\gamma}{\eta}\qquad\qquad
\end{equation*}
\begin{equation*}
    +\frac{1}{\alpha^2(\alpha-1)\Gamma(-\alpha)}\frac{1}{\eta^2}+\frac{\gamma\nu_1}{\eta^2}+(2\mu_1\nu_1+\mu_2\nu_3+\mu_3\nu_2)\frac{1}{\eta^3}+(\nu_1^2+\nu_2\nu_3)\frac{1}{\eta^4}\quad
\end{equation*}
\begin{equation*}
    \le-\frac{\mu_1\gamma}{\eta}+\Big[\frac{1}{\alpha^2(\alpha-1)\Gamma(-\alpha)}+\gamma\nu_1+2\mu_1\nu_1+\mu_2\nu_3+\mu_3\nu_2+\nu_1^2+\nu_2\nu_3\Big]\frac{1}{\eta^2}
\end{equation*}
for $\eta\ge\frac{1}{|\cos{\theta}|}\ge1$. For the last inequality, we use $\frac{1}{\eta^4}\le\frac{1}{\eta^2}$ and $\frac{1}{\eta^3}\le\frac{1}{\eta^2}$ for $\eta\ge1$. Hence
\begin{equation*}
    \rho(\eta)=-\frac{\mu_1\gamma}{\eta}\left\{1-\frac{1}{\mu_1\gamma}\left(\frac{1}{\alpha^2(\alpha-1)\Gamma(-\alpha)}+\gamma\nu_1+2\mu_1\nu_1+\mu_2\nu_3+\mu_3\nu_2+\nu_1^2+\nu_2\nu_3\right)\frac{1}{\eta}\right\}
\end{equation*}
for $\eta\ge\frac{1}{|\cos{\theta}|}$. Thus the proof of Lemma 2.2 is complete.

{\bf Remark 2.1}. In other cases of $\gamma$, there is no information about an upper bound of the largest zero $\eta_N$.

\subsection{Mild solution of FVP \eqref{Eq1}-\eqref{Eq2}}

In this subsection, we will define mild solutions.  According to the definition of the operator $A$, the equality $Ae_n=\lambda_n e_n$ is valid. Hence, in view of the Fourier expansion $u(t)=\sum_{n=1}^{\infty}u_n(t)e_n$, where $u_n(t)=(u(t),e_n),$ Eq. \eqref{Eq1} can be rewritten as
\begin{equation*}
    \Big(\partial_t^{\alpha}\sum\limits_{n=1}^{\infty}u_n(t)e_n,e_n\Big)+\Big(\sum\limits_{n=1}^{\infty}\lambda_nu_n(t)e_n,e_n\Big)=\left(F(t),e_n\right),\quad t\in(0,T).
\end{equation*}
This is equivalent to the equation
\begin{equation*}
    \partial_t^{\alpha}u_n(t)+\lambda_nu_n(t)=F_n(t),
\end{equation*}
where $F_n(t)=\left(F(t),e_n\right).$ By applying the method of solutions of FVPs for fractional ordinary equations in Subsection \eqref{subsection:2.2}, and using the final value data \eqref{Eq2}, we arrive
\begin{equation*}
    u_n(t)=\frac{E_{\alpha,1}(-\lambda_n t^{\alpha})}{\rho(\lambda_n T^{\alpha})}\Big[E_{\alpha,1}(-\lambda_n T^{\alpha})\varphi_n-TE_{\alpha,2}(-\lambda_n T^{\alpha})\psi_n\qquad\qquad\qquad\qquad\qquad\qquad\qquad\qquad
    \end{equation*}
    \begin{equation*}
-E_{\alpha,1}(-\lambda_n T^{\alpha})F_n(r)\star r^{\alpha-1}E_{\alpha,\alpha}(-\lambda_n r^{\alpha})\Big|_{r=T}\qquad\qquad\qquad\qquad\qquad\qquad\quad
  \end{equation*}
    \begin{equation*}
+TE_{\alpha,2}(-\lambda_n T^{\alpha})F_n(r)\star r^{\alpha-2}E_{\alpha,\alpha-1}(-\lambda_n r^{\alpha})\Big|_{r=T}\Big]\qquad\qquad\qquad\qquad\qquad
\end{equation*}
\begin{equation*}
    +\frac{tE_{\alpha,2}(-\lambda_n t^{\alpha})}{\rho(\lambda_n T^{\alpha})}\Big[(\gamma+E_{\alpha,1}(-\lambda_n T^{\alpha}))\psi_n+\lambda_nT^{\alpha-1}E_{\alpha,\alpha}(-\lambda_n T^{\alpha})\varphi_n\quad\qquad
\end{equation*}
\begin{equation*}
-(\gamma+E_{\alpha,1}(-\lambda_n T^{\alpha}))F_n(r)\star r^{\alpha-2}E_{\alpha,\alpha-1}(-\lambda_n r^{\alpha})\Big|_{r=T}\qquad\qquad\qquad\qquad
\end{equation*}
\begin{equation*}
-\lambda_nT^{\alpha-1}E_{\alpha,\alpha}(-\lambda_nT^{\alpha})F_n(r)\star r^{\alpha-1}E_{\alpha,\alpha}(-\lambda_n r^{\alpha})\Big|_{r=T}\Big]\qquad\qquad\qquad\qquad    
\end{equation*}
\begin{equation}\label{Eq21}
+F_n(t)\star t^{\alpha-1}E_{\alpha,\alpha}(-\lambda_n t^{\alpha}),\qquad\qquad\qquad\qquad\qquad\qquad\qquad\qquad\qquad\qquad
\end{equation}
where $\varphi_n=(\varphi,e_n)$ and $\psi_n=(\psi,e_n)$. Thus, we obtain a spectral representation for $u$ as follows:
\begin{equation*}
    u(t)\qquad\qquad\qquad\qquad\qquad\qquad\qquad\qquad\qquad\qquad\qquad\qquad\qquad\qquad\qquad\qquad\qquad\qquad\qquad\qquad
\end{equation*}
\begin{equation*}
    =\sum\limits_{n=1}^{\infty}\Big[E_{\alpha,1}(-\lambda_n T^{\alpha})E_{\alpha,1}(-\lambda_n t^{\alpha})+\lambda_nT^{\alpha-1}E_{\alpha,\alpha}(-\lambda_n T^{\alpha})tE_{\alpha,2}(-\lambda_n t^{\alpha})\Big]\frac{\varphi_n}{\rho(\lambda_nT^{\alpha})}e_n
\end{equation*}
\begin{equation*}
    +\sum\limits_{n=1}^{\infty}\Big[-TE_{\alpha,2}(-\lambda_n T^{\alpha})E_{\alpha,1}(-\lambda_n t^{\alpha})+(\gamma+E_{\alpha,1}(-\lambda_nT^{\alpha}))tE_{\alpha,2}(-\lambda_n t^{\alpha})\Big]\frac{\psi_n}{\rho(\lambda_nT^{\alpha})}e_n
\end{equation*}
\begin{equation*}
    -\sum\limits_{n=1}^{\infty}\Big[E_{\alpha,1}(-\lambda_n T^{\alpha})E_{\alpha,1}(-\lambda_n t^{\alpha})+\lambda_nT^{\alpha-1}E_{\alpha,\alpha}(-\lambda_n T^{\alpha})tE_{\alpha,2}(-\lambda_n t^{\alpha})\Big]\qquad\qquad\qquad
\end{equation*}
\begin{equation*}
    \times\frac{F_n(r)\star r^{\alpha-1}E_{\alpha,\alpha}(-\lambda_nr^{\alpha})\Big|_{r=T}}{\rho(\lambda_nT^{\alpha})}e_n
\end{equation*}
\begin{equation*}
    +\sum\limits_{n=1}^{\infty}\Big[TE_{\alpha,2}(-\lambda_n T^{\alpha})E_{\alpha,1}(-\lambda_n t^{\alpha})-(\gamma+E_{\alpha,1}(-\lambda_nT^{\alpha}))tE_{\alpha,2}(-\lambda_n t^{\alpha})\Big]\qquad\qquad\qquad
\end{equation*}
\begin{equation*}
    \times\frac{F_n(r)\star r^{\alpha-2}E_{\alpha,\alpha-1}(-\lambda_nr^{\alpha})\Big|_{r=T}}{\rho(\lambda_nT^{\alpha})}e_n
\end{equation*}
\begin{equation}\label{Eq22}
    +\sum\limits_{n=1}^{\infty}F_n(t)\star t^{\alpha-1}E_{\alpha,\alpha}(-\lambda_n t^{\alpha})e_n.\qquad\qquad\qquad\qquad\qquad\qquad\qquad\qquad\qquad
\end{equation}
For $g\in C([0,T];H)$ and $\varphi,\psi\in H$, let us denote by $\mathcal{P}_j,\,1\le j\le 5$, the following operators
\begin{equation*}
\begin{cases}
    (\mathcal{P}_1g)(t):=\sum\limits_{n=1}^{\infty}g_n(t)\star t^{\alpha-1}E_{\alpha,\alpha}(-\lambda_n t^{\alpha})e_n,\\
\mathcal{P}_2(t)\varphi:=\sum\limits_{n=1}^{\infty}\Big[E_{\alpha,1}(-\lambda_n T^{\alpha})E_{\alpha,1}(-\lambda_n t^{\alpha})+\lambda_nT^{\alpha-1}E_{\alpha,\alpha}(-\lambda_n T^{\alpha})tE_{\alpha,2}(-\lambda_n t^{\alpha})\Big]\frac{\varphi_n}{\rho(\lambda_nT^{\alpha})}e_n,\\
     \mathcal{P}_3(t)\psi:=\sum\limits_{n=1}^{\infty}\Big[-TE_{\alpha,2}(-\lambda_n T^{\alpha})E_{\alpha,1}(-\lambda_n t^{\alpha})+(\gamma+E_{\alpha,1}(-\lambda_nT^{\alpha}))tE_{\alpha,2}(-\lambda_n t^{\alpha})\Big]\frac{\psi_n}{\rho(\lambda_nT^{\alpha})}e_n,\\
     (\mathcal{P}_4g)(t):=-\mathcal{P}_2(t)(\mathcal{P}_1g)(T),\qquad (\mathcal{P}_5g)(t):=-\mathcal{P}_3(t)(\mathcal{P}_1g)(T).
\end{cases}
\end{equation*}
Then, the solution $u$ can be represented as
\begin{equation}\label{Eq23}
    u(t)=(\mathcal{P}_1F)(t)+\mathcal{P}_2(t)\varphi+\mathcal{P}_3(t)\psi+(\mathcal{P}_4F)(t)+(\mathcal{P}_5F)(t),
\end{equation}
for $t\in(0,T).$

{\bf Definition 2.1.} \textit{For any given $\alpha\in(1,2),\gamma\not=0,\varphi,\psi\in H$ and $F\in C([0,T],H)$, a function $u\in C([0,T];H)\cap C((0,T];D(A))$ is called a mild solution of \eqref{Eq1}, \eqref{Eq2} if it satisfies the integral equation \eqref{Eq23}.} 

We define the operator valued function $Y_A(t)$ by
\begin{equation*}
    Y_A(t)h=\sum\limits_{n=1}^{\infty}(h,e_n)t^{\alpha-1}E_{\alpha,\alpha}(-\lambda_nt^{\alpha})e_n
\end{equation*}
for $h\in H$. Moreover,  we set
\begin{equation*}
    \Lambda=\Lambda(\alpha,A):=\bigcup\limits_{n=1}^{\infty}\left\{\Big(\frac{\eta_0}{\lambda_n}\Big)^{\frac{1}{\alpha}},...,\Big(\frac{\eta_N}{\lambda_n}\Big)^{\frac{1}{\alpha}}\right\}.
\end{equation*}
We note that $\Lambda$ is a countably infinite set. Since $\lim\limits_{n\rightarrow\infty}\lambda_n=\infty$, the set $\Lambda$ has an accumulation point $0$, but we can readily verify
\begin{equation*}
\Lambda\subset\left[0,\Big(\frac{\eta_N}{\lambda_1}\Big)^{\frac{1}{\alpha}}\right].
\end{equation*}

Therefore, we have 

{\bf Lemma 2.3.} \textit{Let $\gamma\not=0$ and $T\not\in\Lambda$. If $\varphi,\psi\in D(A^{1-\frac{1}{\alpha}})$ and $F\in L^{\infty}(0,T;D(A^{1-\frac{1}{\alpha}}))$, then there exists a constant $C_0>0$ such that
\begin{equation}\label{lemma2.3}
    \|u\|_{C([0,T];H)}\le C_0\left(\|\varphi\|_{D(A^{1-\frac{1}{\alpha}})}+\|\psi\|_{D(A^{1-\frac{1}{\alpha}})}+\|F\|_{L^{\infty}(0,T;D(A^{1-\frac{1}{\alpha}})}\right),
\end{equation}
where $C_0$ depends only $\gamma,\alpha,\Omega$ and $T$.
}

{\bf Proof.} By Lemma 1.1 and the definition $(\mathcal{P}_1F)(t)$, also by the generalized Minkovskii inequality, we have that
\begin{equation*}
    \|(\mathcal{P}_1F)(t)\|\le \int_0^t\|Y_A(t-s)\|\|F(s)\|ds\le c\int_0^t(t-s)^{\alpha-1}\|F(s)\|ds
    \le c\frac{t^{\alpha}}{\alpha}\|F\|_{L^{\infty}(0,T;H)}.
\end{equation*}
In addition, by Lemma 1.1 and \eqref{Eq15}, the norm $\|\mathcal{P}_2(t)\varphi\|$ can be estimated as
\begin{equation}\label{Eq24}
    \|\mathcal{P}_2(t)\varphi\|^2\le2c^4c_1^{-2}\sum\limits_{n=1}^{\infty}\left|\frac{\lambda_nT^{\alpha}}{1+\lambda_nT^{\alpha}}\frac{1}{1+\lambda_nt^{\alpha}}\varphi_n\right|^2+2c^4c_1^{-2}\sum\limits_{n=1}^{\infty}\left|\frac{\lambda_nT^{2\alpha-1}}{1+\lambda_nT^{\alpha}}\frac{(\lambda_nt^{\alpha})^{\frac{1}{\alpha}}}{1+\lambda_nt^{\alpha}}\lambda_n^{1-\frac{1}{\alpha}}\varphi_n\right|^2
\end{equation}
for all $t\in[0,T]$. Applying
\begin{equation}\label{Eq25}
    \sup\limits_{0< y<\infty}\frac{y^{\theta}}{1+y}=\frac{\left(\frac{\theta}{1-\theta}\right)^{\theta}}{1+\frac{\theta}{1-\theta}},\qquad0<\theta<1
\end{equation}
in \eqref{Eq24}, we have
\begin{equation}\label{Eq26}
    \|\mathcal{P}_2(t)\varphi\|\le c_2\|\varphi\|_{D(A^{1-\frac{1}{\alpha}})},\quad t\in[0,T],
\end{equation}
here we have used from the $D(A^{1-\frac{1}{\alpha}})\subset H$ for $1<\alpha<2$. Similarly to \eqref{Eq26}, for $\mathcal{P}_3(t)\psi$, we have
\begin{equation}\label{Eq27}
       \|\mathcal{P}_3(t)\psi\|\le c_3\|\psi\|_{D(A^{1-\frac{1}{\alpha}})},\quad t\in[0,T].
\end{equation}
Now, we proceed to estimate $\|(\mathcal{P}_4F)(t)\|$ by using the same technique as in \eqref{Eq24}. As a consequence of
\begin{equation*}
     (\mathcal{P}_4F)(t):=-\mathcal{P}_2(t)(\mathcal{P}_1F)(T),
\end{equation*}
we can obtain the following estimates
\begin{equation*}
    \|(\mathcal{P}_4F)(t)\|^2\le2c^6c_1^2 \sum\limits_{n=1}^{\infty}\left|\frac{\lambda_nT^{\alpha}}{1+\lambda_nT^{\alpha}}\frac{1}{1+\lambda_nt^{\alpha}}\int_0^TF_n(s)\frac{(T-s)^{\alpha-1}}{1+\lambda_n(T-s)^{\alpha}}ds\right|^2\qquad\qquad\qquad\qquad\qquad\qquad
\end{equation*}
\vspace{0.2cm}
\begin{equation*}
    +2c^6c_1^2 T^{2(\alpha-1)}\sum\limits_{n=1}^{\infty}\left|\frac{\lambda_nT^{\alpha}}{1+\lambda_nT^{\alpha}}\frac{(\lambda_nt^{\alpha})^{\frac{1}{\alpha}}}{1+\lambda_nt^{\alpha}}\int_0^TF_n(s)\frac{(\lambda_n(T-s)^{\alpha})^{\frac{\alpha-1}{\alpha}}}{1+\lambda_n(T-s)^{\alpha}}ds\right|^2
\end{equation*}
\vspace{0.2cm}
\begin{equation*}
    \stackrel{{g.M.ineq.}}{\le}2c^6c_1^2\left|\int_0^T\Big(\sum\limits_{n=1}^{\infty}F_n^2(s)\Big)^{\frac{1}{2}}(T-s)^{\alpha-1}ds\right|^2+\frac{2c^6c_1^2}{\alpha^4}(\alpha-1)^{\frac{4\alpha-4}{\alpha}}T^{2(\alpha-1)}\left|\int_0^T\Big(\sum\limits_{n=1}^{\infty}F_n^2(s)\Big)^{\frac{1}{2}}ds\right|^2
\end{equation*}
\vspace{0.2cm}
\begin{equation}\label{Eq28}
    \le c_4^2\|F\|^2_{L^{\infty}(0,T;H)}.
\end{equation}
A simple computation shows that
\begin{equation*}
    \|(\mathcal{P}_5F)(t)\|^2\le2c^6c_1^2 \sum\limits_{n=1}^{\infty}\left|\frac{\lambda_nT^{\alpha+1}}{1+\lambda_nT^{\alpha}}\frac{1}{1+\lambda_nt^{\alpha}}\int_0^TF_n(s)(T-s)^{\alpha-2}ds\right|^2\qquad\qquad\qquad\qquad
\end{equation*}
\vspace{0.2cm}
\begin{equation*}
    +2c^6c_1^2 T^{2\alpha}\sum\limits_{n=1}^{\infty}\left|\left(\gamma+\frac{1}{1+\lambda_nT^{\alpha}}\right)\frac{(\lambda_nt^{\alpha})^{\frac{1}{\alpha}}}{1+\lambda_nt^{\alpha}}\int_0^T\lambda_n^{1-\frac{1}{\alpha}}F_n(s)(T-s)^{\alpha-2}ds\right|^2
\end{equation*}
\vspace{0.2cm}
\begin{equation}\label{Eq29}
    \le c_5^2\|F\|^2_{L^{\infty}(0,T;D(A^{1-\frac{1}{\alpha}})}.
\end{equation}

Thus the proof of Lemma 2.3 is complete.

Based on Lemma 2.3, we consider solutions existence, uniqueness and regularity in the following theorem. 

{\bf Theorem 2.1.} \textit{Let $\gamma\not=0$ and $T\not\in\Lambda$. 1) If $\varphi,\psi\in D(A^{1-\frac{1}{\alpha}})$ and $F\in L^{\infty}(0,T;D(A^{1-\frac{1}{\alpha}}))$, then FVP \eqref{Eq1}-\eqref{Eq2} has a unique solution $u\in C([0,T];H)$. Furthermore, there exists a constant $C_0>0$ such that
\begin{equation}\label{Th1.1}
    \|u\|_{C([0,T];H)}\le C_0\left(\|\varphi\|_{D(A^{1-\frac{1}{\alpha}})}+\|\psi\|_{D(A^{1-\frac{1}{\alpha}})}+\|F\|_{L^{\infty}(0,T;D(A^{1-\frac{1}{\alpha}})}\right),
\end{equation}
where $C_0$ depends only $\gamma,\alpha,\Omega$ and $T$.\\
2) If $\varphi,\psi\in D(A)$ and $F\in L^{\infty}(0,T;D(A))$, then FVP \eqref{Eq1}-\eqref{Eq2} has a unique solution $u\in C([0,T];H)\cap C((0,T];D(A))\cap  C^1([0,T];D(A))$ to \eqref{Eq1}-\eqref{Eq2} such that $\partial_t^{\alpha}u\in C((0,T];H)$. Moreover, there exists a constant $C_1>0$ such that
\begin{equation*}
    \|u(t)\|_{D(A)}+\|\partial_t^{\alpha}u(t,\cdot
    )\|\le C_1\Big[\left(t^{-\alpha}+t^{1-\alpha}\right)(\|\varphi\|_{D(A)}+\|\psi\|_{D(A)})\qquad\qquad\qquad
\end{equation*}
\begin{equation}\label{Th1.2}
\qquad\qquad\qquad\qquad+\left(1+t^{\alpha}+t^{1-\alpha}\right)\|F\|_{L^{\infty}(0,T;D(A))}\Big],\quad t>0.
\end{equation}\\
3) Moreover, if $\varphi\in D(A)$, $\psi\in D(A^{2-\frac{1}{\alpha}})$ and $F\in AC([0,T];D(A^{2-\frac{1}{\alpha}}))$, then  there is a unique solution to FVP \eqref{Eq1}-\eqref{Eq2} satisfying $u\in AC^2([0,T];H)\cap C([0,T];D(A))$ and $\partial_t^{\alpha}u\in L^{\infty}(0,T;H),$ and we have the following estimate
\begin{equation}\label{Th1.3}
    \|\partial_t^{\alpha}u\|_{L^{\infty}(0,T;H)}\le C_2\left(\|\varphi\|_{D(A)}+\|\psi\|_{D(A^{2-\frac{1}{\alpha}})}+\|F\|_{L^{\infty}(0,T;D(A^{2-\frac{1}{\alpha}}))}\right),
\end{equation}
where $C_2>0$ is depending on $\alpha,\gamma,\Omega$ and $T$.\\
4) Furthermore, if $\varphi,\, \psi\in D(A)$ and $F\in C^{0,\alpha-1}([0,T];D(A))$, then for the solution $u$ given by \eqref{Eq23}, we have:
 for every $\delta>0$,
\begin{equation*}
     \|Au\|_{C^{\alpha-1}([0,T];H)}+\|\partial_t^{\alpha}u\|_{C^{\alpha-1}([0,T];H)}\le C_3\Big[\|F\|_{C^{0,\alpha-1}([0,T];D(A))}\qquad\qquad\qquad\qquad
\end{equation*}
\begin{equation}\label{Th1.5}
  \qquad\qquad\qquad\qquad+\Big(\frac{1}{\delta}+\frac{1}{\delta^{\alpha-1}}\Big)(\|\varphi\|_{D(A)}+\|\psi\|_{D(A)}+\|F\|_{C([0,T];D(A))}\Big];
\end{equation}
where $C_3$ is positive constant, doesn't depend on $\varphi,\psi$ and $F$.
}

{\bf Proof.}  The proof of Part 1) is the same as Lemma 2.3. The uniqueness of $u$ depends on the uniqueness of the  \eqref{Eq5}, \eqref{Eq6}. Indeed, 
suppose we have two solutions:  $u_1(t),u_2(t)$ and set $v(t)=u_1(t)-u_2(t)$. Then, we have
\begin{equation}\label{un1}
    \begin{cases}
        \partial_{t}^{\alpha}v(t)+Av(t)=0,&0<t<T,\\
         \gamma v(0)+v(T)=0,\quad \partial_tv(T)=0.
    \end{cases}
\end{equation}
Set $v_n(t)=(v(t),e_n)$. It follows from Eq. in \eqref{un1} that for any $n\in\mathbf{N}$
\begin{equation*}
    \partial_t^{\alpha}v_n(t)=(\partial_t^{\alpha}v(t,x),e_n)=-(Av(t,x),e_n)=-(v(t,x),Ae_n)=-\lambda_nv_n(t)
\end{equation*}
and the final time conditions imply
\begin{equation*}
    \gamma v_n(0)+v_n(T)=0,\quad \partial_tv_n(T)=0.
\end{equation*}
Therefore, we get the following FVP for $v_n(t)$
\begin{equation}\label{un2}
    \begin{cases}
        \partial_{t}^{\alpha}v_n(t)+\lambda_nv_n(t)=0,&0<t<T,\\
         \gamma v_n(0)+v_n(T)=0,\quad \partial_tv_n(T)=0.
    \end{cases}
\end{equation}
By our assumption, $T$ does not belong to the set $\Lambda$. Then, by \eqref{Eq14}, the problem \eqref{un2} has only zero solution, i.e., $v_n(t)\equiv0$. Then by the Parseval's identity, we obtain $v(t)\equiv0$ for all $[0, T]\time\Omega$. Hence, the uniqueness of the solution is proved.  

Now, we proceed to prove Part 2). By Lemma 1.1, we have
\begin{equation*}
    \|(\mathcal{P}_1F)(t)\|_{D(A)}=\|A(\mathcal{P}_1F)(t)\|\stackrel{{g.M.ineq.}}{\le}\int_0^t\|Y_A(t-s)AF(s)\|ds\qquad\qquad\qquad\qquad\qquad
\end{equation*}
\begin{equation}\label{Eq29.0}
\le c\int_0^t(t-s)^{\alpha-1}\|F(s)\|_{D(A)}ds\le c\frac{t^{\alpha}}{\alpha}\|F\|_{L^{\infty}(0,T;D(A))}.    
\end{equation}
Similarly, from $\|\mathcal{P}_2(t)\varphi\|_{D(A)}=\|A\mathcal{P}_2(t)\varphi\|$ and the same way as in the proof of \eqref{Eq26}, we have
\begin{equation*}
    \|\mathcal{P}_2(t)\varphi\|^2_{D(A)}\qquad\qquad\qquad\qquad\qquad\qquad\qquad\qquad\qquad\qquad\qquad\qquad\qquad\qquad\qquad
\end{equation*}
\begin{equation*}
\le2c^4c_1^{-2}\sum\limits_{n=1}^{\infty}\left|\frac{\lambda_nT^{\alpha}}{1+\lambda_nT^{\alpha}}\frac{\lambda_nt^{\alpha}}{1+\lambda_nt^{\alpha}}t^{-\alpha}\varphi_n\right|^2+2c^4c_1^{-2}\sum\limits_{n=1}^{\infty}\left|\frac{\lambda_nT^{2\alpha-1}}{1+\lambda_nT^{\alpha}}\frac{\lambda_nt^{\alpha}}{1+\lambda_nt^{\alpha}}t^{1-\alpha}\lambda_n\varphi_n\right|^2 
\end{equation*}
\vspace{0.2cm}
\begin{equation}\label{Eq30}
    \le c_6^2(t^{-2\alpha}+t^{2-2\alpha})\|\varphi\|^2_{D(A)},\quad t>0.
\end{equation}
For $\|\mathcal{P}_3(t)\psi\|_{D(A)}$, we have
\begin{equation*}
    \|\mathcal{P}_3(t)\psi\|^2_{D(A)}\qquad\qquad\qquad\qquad\qquad\qquad\qquad\qquad\qquad\qquad\qquad\qquad\qquad\qquad\qquad
\end{equation*}
\begin{equation*}
\le2c^4c_1^{-2}\sum\limits_{n=1}^{\infty}\left|\frac{\lambda_nT^{\alpha+1}}{1+\lambda_nT^{\alpha}}\frac{\lambda_nt^{\alpha}}{1+\lambda_nt^{\alpha}}t^{-\alpha}\psi_n\right|^2+2c^4c_1^{-2}T^{2\alpha}\sum\limits_{n=1}^{\infty}\left|\Big(\gamma+\frac{1}{1+\lambda_nT^{\alpha}}\Big)\frac{\lambda_nt^{\alpha}}{1+\lambda_nt^{\alpha}}t^{1-\alpha}\lambda_n\psi_n\right|^2 
\end{equation*}
\vspace{0.2cm}
\begin{equation}\label{Eq31}
    \le c_6^2(t^{-2\alpha}+t^{2-2\alpha})\|\psi\|^2_{D(A)},\quad t>0.
\end{equation}
Now, we will estimate the norm $\|(\mathcal{P}_4F)(t)\|_{D(A)},$ which will use the same as in \eqref{Eq30}. Indeed, by Lemma 1.1 and \eqref{Eq25}, we see that
\begin{equation*}
    \|(\mathcal{P}_4F)(t)\|_{D(A)}^2\le2c^6c_1^{-2} \sum\limits_{n=1}^{\infty}\left|\frac{\lambda_nT^{\alpha}}{1+\lambda_nT^{\alpha}}\frac{1}{1+\lambda_nt^{\alpha}}\int_0^T\lambda_n^{\frac{1}{\alpha}}F_n(s)\frac{(\lambda_n(T-s)^{\alpha})^{\frac{\alpha-1}{\alpha}}}{1+\lambda_n(T-s)^{\alpha}}ds\right|^2\qquad\qquad\qquad\qquad\qquad\qquad
\end{equation*}
\vspace{0.2cm}
\begin{equation*}
    +2c^6c_1^{-2}T^{2(\alpha-1)}\sum\limits_{n=1}^{\infty}\left|\frac{\lambda_nT^{\alpha}}{1+\lambda_nT^{\alpha}}\frac{(\lambda_nt^{\alpha})^{\frac{1}{\alpha}}}{1+\lambda_nt^{\alpha}}\int_0^T\lambda_nF_n(s)\frac{(\lambda_n(T-s)^{\alpha})^{\frac{\alpha-1}{\alpha}}}{1+\lambda_n(T-s)^{\alpha}}ds\right|^2
\end{equation*}
\vspace{0.2cm}
\begin{equation*}
    \stackrel{{g.M.ineq.}}{\le}2c^6c_1^{-2}\left|\int_0^T\Big(\sum\limits_{n=1}^{\infty}(\lambda_n^{\frac{1}{\alpha}}F_n(s))^2\Big)^{\frac{1}{2}}ds\right|^2+\frac{2c^6c_1^{-2}}{\alpha^4}(\alpha-1)^{\frac{4\alpha-4}{\alpha}}T^{2(\alpha-1)}\left|\int_0^T\Big(\sum\limits_{n=1}^{\infty}(\lambda_nF_n(s))^2\Big)^{\frac{1}{2}}ds\right|^2
\end{equation*}
\vspace{0.2cm}
\begin{equation}\label{Eq32}
    \le c_7^2\|F\|^2_{L^{\infty}(0,T;D(A^{1/\alpha}))}+c_8^2\|F\|^2_{L^{\infty}(0,T;D(A))}.
\end{equation}
Similarly computation shows that
\begin{equation}\label{Eq33}
    \|(\mathcal{P}_5F)(t)\|\le c_9(1+t^{1-\alpha})\|F\|_{L^{\infty}(0,T;D(A))},\quad t>0.
\end{equation}
Bring, \eqref{Eq29.0}-\eqref{Eq33}, we have
\begin{equation*}
     \|u(t)\|_{D(A)}\le c_{10}\left(t^{-\alpha}+t^{1-\alpha}\right)(\|\varphi\|_{D(A)}+\|\psi\|_{D(A)})\qquad\qquad\qquad\qquad
\end{equation*}
\begin{equation}\label{Eq34}
\qquad\qquad+c_{11}\left(1+t^{\alpha}+t^{1-\alpha}\right)\|F\|_{L^{\infty}(0,T;D(A))},\quad t>0.
\end{equation}
Furthermore,  by \eqref{Eq1} and \eqref{Eq29.0}-\eqref{Eq33}, we have
\begin{equation*}
   \|Au(t)\|^2_{H} =\|u\|^2_{D(A)}
\end{equation*}
\begin{equation}\label{Eq35}
\le 
c_{10}^{2'}\left(t^{-2\alpha}+t^{2-2\alpha}\right)(\|\varphi\|^2_{D(A)}+\|\psi\|^2_{D(A)})+c_{11}^{2'}\left(1+t^{2\alpha}+t^{2-2\alpha}\right)\|F\|^2_{L^{\infty}(0,T;D(A))},\quad t>0.
\end{equation}
In addition, from Eq.\eqref{Eq1} and \eqref{Eq35}, we obtain
\begin{equation*}
    \|\partial_t^{\alpha}u(t)\|\le 
c_{12}\left(t^{-\alpha}+t^{1-\alpha}\right)(\|\varphi\|_{D(A)}+\|\psi\|_{D(A)})\qquad\qquad\qquad\qquad
\end{equation*}
\begin{equation}\label{Eq36}
\qquad\qquad+c_{13}\left(1+t^{\alpha}+t^{1-\alpha}\right)\|F\|_{L^{\infty}(0,T;D(A))},\quad t>0.
\end{equation}

Now, we will show that the $u$ in $ C^1([0,T];D(A))$. Indeed, by Lemma 3, we have
\begin{equation*}
    \partial_t(\mathcal{P}_1F)(t)=\sum\limits_{n=1}^{\infty}F_n(t)\star t^{\alpha-2}E_{\alpha,\alpha-1}(t,\lambda_n)e_n,
\end{equation*}
then, by the generalized Minkowski inequality, we get
\begin{equation*}
    \|\partial_t(\mathcal{P}_1F)(t)\|\le\int_0^t\left\|\sum\limits_{n=1}^{\infty}(t-s)^{\alpha-2}E_{\alpha,\alpha-1}(-\lambda_n(t-s)^{\alpha})F_n(s)\right\|ds
\end{equation*}
\begin{equation*}
\quad\qquad\le\int_0^t\left[\sum\limits_{n=1}^{\infty}\Big(F_n(s)\frac{(t-s)^{\alpha-2}}{1+\lambda_n(t-s)^{\alpha}}\Big)^2\right]^{1/2}ds
\end{equation*}
\begin{equation}\label{Eq37}
\le \frac{t^{\alpha-1}}{\alpha-1}\|F\|_{L^{\infty}(0,T;H)}.\qquad\qquad
\end{equation}
Again using Lemma 1.3, we have
\begin{equation*}
    \partial_t(\mathcal{P}_2(t)\varphi)=\sum\limits_{n=1}^{\infty}\Big[-\lambda_nt^{\alpha-1}E_{\alpha,\alpha}(-\lambda_nt^{\alpha})E_{\alpha,1}(-\lambda_nT^{\alpha})\qquad\qquad\qquad
\end{equation*}
\begin{equation*}
\qquad\qquad\qquad\qquad\qquad+\lambda_nT^{\alpha-1}E_{\alpha,\alpha}(-\lambda_nT^{\alpha})E_{\alpha,1}(-\lambda_nt^{\alpha})\Big]\frac{\varphi_n}{\rho(\lambda_nT^{\alpha})}e_n,
\end{equation*}
then, by virtue of \eqref{Eq15} and \eqref{Eq25}, we obtain
\begin{equation*}
    \|\partial_t(\mathcal{P}_2(t)\varphi)\|^2=
\sum\limits_{n=1}^{\infty}\bigg|\Big[-\lambda_nt^{\alpha-1}E_{\alpha,\alpha}(-\lambda_nt^{\alpha})E_{\alpha,1}(-\lambda_nT^{\alpha})\qquad\qquad\qquad\qquad\qquad\qquad\qquad
\end{equation*}
\begin{equation*}
+\lambda_nT^{\alpha-1}E_{\alpha,\alpha}(-\lambda_nT^{\alpha})E_{\alpha,1}(-\lambda_nt^{\alpha})\Big]\frac{\varphi_n}{\rho(\lambda_nT^{\alpha})}\bigg|^2\le 2c^4c_1^{-2}\sum\limits_{n=1}^{\infty}\frac{\lambda_n^2T^{2\alpha}}{(1+\lambda_nT^{\alpha})^2}\left(\frac{(\lambda_nt^{\alpha})^{\frac{\alpha-1}{\alpha}}}{1+\lambda_nt^{\alpha}}\right)^2\lambda_n^{\frac{2}{\alpha}}\varphi_n^2
\end{equation*}
\begin{equation*}
    +2c^4c_1^{-2}\sum\limits_{n=1}^{\infty}\frac{\lambda_n^2}{(1+\lambda_nt^{\alpha})^2}\left(\frac{\lambda_n^2T^{2\alpha-1}}{1+\lambda_nT^{\alpha}}\right)^2\varphi_n^2.
\end{equation*}
By \eqref{Eq19} and $D(A)\subset D(A^{\frac{1}{\alpha}})$ for $1<\alpha<2$, we have
\begin{equation}\label{Eq38}
    \|\partial_t(\mathcal{P}_2(t)\varphi)\|\le c_{14}\|\varphi\|_{D(A)}.
\end{equation}
Similarly, we have
\begin{equation}\label{Eq39}
    \|\partial_t(\mathcal{P}_3(t)\psi)\|\le c_{15}\|\psi\|_{D(A)}.
\end{equation}
Furthermore, by defining $\mathcal{P}_4F(t)$ and using Lemma 3, we have
\begin{equation*}
    \partial_t(\mathcal{P}_4F)(t)=-\partial_t(\mathcal{P}_2(t)F(T))\qquad\qquad\qquad\qquad\qquad\qquad\qquad\qquad\qquad\qquad\qquad\qquad\qquad\qquad
\end{equation*}
\begin{equation*}
    =\sum\limits_{n=1}^{\infty}\frac{\lambda_n}{\rho(\lambda_nT^{\alpha})}\Big[E_{\alpha,1}(-\lambda_n T^{\alpha})t^{\alpha-1}E_{\alpha,\alpha}(-\lambda_n t^{\alpha})-T^{\alpha-1}E_{\alpha,\alpha}(-\lambda_n T^{\alpha})E_{\alpha,1}(-\lambda_n t^{\alpha})\Big]
\end{equation*}
\begin{equation*}
    \times F_n(r)\star r^{\alpha-1}E_{\alpha,\alpha}(-\lambda_nr^{\alpha})\Big|_{r=T}e_n
\end{equation*}
\begin{equation*}
    :=\mathcal{I}_{11}F(t)+\mathcal{I}_{12}F(t).\qquad\quad\quad\quad
\end{equation*}
Then, similarly to \eqref{Eq31}, we get
\begin{equation*}
    \|\mathcal{I}_{11}F(t)\|^2=\sum\limits_{n=1}^{\infty}\left|\frac{\lambda_n}{\rho(\lambda_nT^{\alpha})}E_{\alpha,1}(-\lambda_nT^{\alpha})t^{\alpha-1}E_{\alpha,\alpha}(-\lambda_nt^{\alpha})\int_0^TF_n(s)(T-s)^{\alpha-1}E_{\alpha,\alpha}(-\lambda_n(T-s)^{\alpha})ds\right|^2
\end{equation*}
\begin{equation*}
    \le c^6c_1^{-2}\sum\limits_{n=1}^{\infty}\frac{\lambda_n^{2}T^{2\alpha}}{(1+\lambda_nT^{\alpha})^2}\left(\frac{(\lambda_nt^{\alpha})^{\frac{\alpha-1}{\alpha}}}{1+\lambda_nt^{\alpha}}\right)^2\left|\int_0^TF_n(s)\frac{(\lambda_n(T-s)^{\alpha})^{\frac{1}{\alpha}}}{1+\lambda_n(T-s)^{\alpha}}(T-s)^{\alpha-2}ds\right|^2
\end{equation*}
\begin{equation*}
    \le c^6c_1^{-2}\left(\alpha-1\right)^{\frac{4\alpha-4}{\alpha}}\left(\frac{1}{\alpha}\right)^{4}\left|\int_0^T\Big(\sum\limits_{n=1}^{\infty}\Big(F(s,\cdot),e_n\Big)^2\Big)^{1/2}(T-s)^{\alpha-2}ds\right|^2
\end{equation*}
\begin{equation}\label{Eq40}
    \le c_{16}^2\|F\|^2_{L^{\infty}(0,T;H)}.
\end{equation}
Similarly, we have
\begin{equation}\label{Eq41}
    \|\mathcal{I}_{12}F(t,\cdot)\|\le c_{17}\|F\|_{L^{\infty}(0,T;D(A^{\frac{1}{\alpha}}))}.
\end{equation}
Using Lemma 1.3, we get
\begin{equation*}
    \partial_t(\mathcal{P}_5F)(t)=-\partial_t(\mathcal{P}_3(t)F(T))\qquad\qquad\qquad\qquad\qquad\qquad\qquad\qquad\qquad\qquad\qquad\qquad\qquad\qquad
\end{equation*}
\begin{equation*}
    =\sum\limits_{n=1}^{\infty}\frac{\lambda_n}{\rho(\lambda_nT^{\alpha})}\Big[TE_{\alpha,2}(-\lambda_n T^{\alpha})t^{\alpha-1}E_{\alpha,\alpha}(-\lambda_n t^{\alpha})+(\gamma+E_{\alpha,1}(-\lambda_n T^{\alpha}))E_{\alpha,1}(-\lambda_n t^{\alpha})\Big]
\end{equation*}
\begin{equation*}
    \times F_n(r)\star r^{\alpha-2}E_{\alpha,\alpha-1}(-\lambda_nr^{\alpha})\Big|_{r=T}e_n
\end{equation*}
\begin{equation*}
    :=\mathcal{I}_{13}F(t)+\mathcal{I}_{14}F(t).\qquad\quad\quad\quad
\end{equation*}
Then, similarly to \eqref{Eq33}, we get
\begin{equation*}
    \|\mathcal{I}_{13}F(t)\|^2\qquad\qquad\qquad\qquad\qquad\qquad\qquad\qquad\qquad\qquad\qquad\qquad\qquad\qquad\qquad\qquad\qquad\qquad\qquad
\end{equation*}
\begin{equation*}
    =\sum\limits_{n=1}^{\infty}\left|\frac{\lambda_nT}{\rho(\lambda_nT^{\alpha})}E_{\alpha,2}(-\lambda_nT^{\alpha})t^{\alpha-1}E_{\alpha,\alpha}(-\lambda_nt^{\alpha})\int_0^TF_n(s)(T-s)^{\alpha-2}E_{\alpha,\alpha-1}(-\lambda_n(T-s)^{\alpha})ds\right|^2
\end{equation*}
\begin{equation*}
    \le c^6c_1^{-2}\sum\limits_{n=1}^{\infty}\frac{\lambda_n^{4}T^{2\alpha+2}}{(1+\lambda_nT^{\alpha})^2}\frac{t^{2(\alpha-1)}}{(1+\lambda_nt^{\alpha})^2}\left|\int_0^TF_n(s)\frac{(T-s)^{\alpha-2}}{1+\lambda_n(T-s)^{\alpha}}ds\right|^2
\end{equation*}
\begin{equation*}
    \le c_0^6c_1^2T^{2}\left(\alpha-1\right)^{\frac{2\alpha-2}{\alpha}}\left(\frac{1}{\alpha}\right)^{2}\left|\int_0^T\Big(\sum\limits_{n=1}^{\infty}\Big(\lambda_n^{\frac{1}{\alpha}}F_n(s)\Big)^2\Big)^{1/2}(T-s)^{\alpha-2}ds\right|^2
\end{equation*}
\begin{equation}\label{Eq42}
    \le c_{18}^2\|F\|^2_{L^{\infty}(0,T;D(A^{\frac{1}{\alpha}}))}.
\end{equation}
For $\mathcal{I}_{14}F(t)$ similar to \eqref{Eq42}, we obtain
\begin{equation}\label{Eq43}
    \|\mathcal{I}_{14}F(t)\|\le c_{19}\|F\|_{L^{\infty}(0,T;D(A))}.
\end{equation}

In the following, we improve the regularity of the solution $u$ such that the $\partial_t^{\alpha}u$ is well-defined by the Caputo fractional derivative. By Lemma 1.2, we have $(\mathcal{P}_1F)(t)\in AC([0,T];H)$ and
\begin{equation*}
    \partial_t^{\alpha}\left(\mathcal{P}_1F(t)\right)=\sum\limits_{n=1}^{\infty}\Big(F_n(t)-\lambda_n\int_0^tF_n(s)(t-s)^{\alpha-1}E_{\alpha,\alpha}(-\lambda_n(t-s)^{\alpha})ds\Big)e_n
\end{equation*}
\begin{equation}\label{Eq44}
    =F(t)-A(\mathcal{P}_1F)(t).\qquad\qquad\qquad\qquad\qquad\qquad\quad
\end{equation}
We see from the last formula and \eqref{Eq29.0} that $\partial_t^{\alpha}(\mathcal{P}_1F)\in C([0,T];H)$. From Lemma 1.3, we know the second derivative for $\mathcal{P}_2(t)\varphi$ satisfy
\begin{equation*}
    \partial_{tt}(\mathcal{P}_2(t)\varphi)=\sum\limits_{n=1}^{\infty}\Big[-\lambda_nt^{\alpha-2}E_{\alpha,\alpha-1}(-\lambda_nt^{\alpha})E_{\alpha,1}(-\lambda_nT^{\alpha})\qquad\qquad\qquad\qquad
\end{equation*}
\begin{equation*}
\qquad\qquad\qquad\qquad\qquad-\lambda_n^2T^{\alpha-1}E_{\alpha,\alpha}(-\lambda_nT^{\alpha})t^{\alpha-1}E_{\alpha,\alpha}(-\lambda_nt^{\alpha})\Big]\frac{\varphi_n}{\rho(\lambda_nT^{\alpha})}e_n,
\end{equation*}
then by \eqref{Eq15} and \eqref{Eq19}, we have
\begin{equation*}
\|\partial_{tt}\left(\mathcal{P}_2(t)\varphi\right)\|^2\le 2c^4c_1^{-2}\sum\limits_{n=1}^{\infty}\left|\frac{t^{\alpha-2}}{1+\lambda_nt^{\alpha}}\frac{1}{1+\lambda_nT^{\alpha}}\lambda_n^2\varphi_n\right|^2 +2c^4c_1^{-2}\sum\limits_{n=1}^{\infty}\left|\frac{t^{\alpha-1}}{1+\lambda_nt^{\alpha}}\frac{\lambda_nT^{\alpha-1}}{1+\lambda_nT^{\alpha}}\lambda_n^2\varphi_n\right|^2 
\end{equation*}
\begin{equation*}
    \le 2c^4c_1^{-2}T^{-2\alpha}t^{2(\alpha-2)}\sum\limits_{n=1}^{\infty}\left|\lambda_n\varphi_n\right|+2c_1^2c_0^4T^{-2}t^{\alpha-2}\sum\limits_{n=1}^{\infty}\left|\lambda_n\varphi_n\right|^2\left(\frac{(\lambda_nt^{\alpha})^{1/2}}{1+\lambda_nt^{\alpha}}\right)^2
\end{equation*}
\begin{equation}\label{Eq45}
    \le c_{20}^2t^{2(\alpha-2)}\|\varphi\|^2_{D(A)}+c_{21}^2t^{\alpha-2}\|\varphi\|^2_{D(A)},
\end{equation}
thus, note that $1<\alpha<2$, we have $\partial_{tt}\left(\mathcal{P}_2(t)\varphi\right)\in L^1(0,T;H)$ and $(\mathcal{P}_2(t)\varphi)\in AC^2([0,T];H)$.

From Lemma 1.3, we know the second derivative for $\mathcal{P}_3(t)\psi$ satisfy
\begin{equation*}
\partial_{tt}\left(\mathcal{P}_3(t)\psi\right)=\sum\limits_{n=1}^{\infty}\lambda_n\Big[TE_{\alpha,2}(-\lambda_nT^{\alpha})t^{\alpha-2}E_{\alpha,\alpha-1}(-\lambda_nt^{\alpha})\qquad\qquad\qquad\qquad\qquad
\end{equation*}
\begin{equation*}
   \qquad\qquad\qquad -\left(\gamma+E_{\alpha,1}(-\lambda_nT^{\alpha})\right)t^{\alpha-1}E_{\alpha,\alpha}(-\lambda_nt^{\alpha})\Big]\frac{\psi_n}{\rho(\lambda_nT^{\alpha})}e_n,
\end{equation*}
then by Lemma 1.1 and \eqref{Eq25}, we have
\begin{equation}\label{Eq46}
\|\partial_{tt}\left(\mathcal{P}_3(t)\psi\right)\|\le c_{22}t^{\alpha-2}\|\psi\|_{D(A)}+c_{20}t^{\alpha-2}\|\psi\|_{D(A^{2-\frac{1}{\alpha}})},
\end{equation}
thus, note that $1<\alpha<2$, we have $\partial_{tt}\left(\mathcal{P}_3(t)\psi\right)\in L^1(0,T;H)$ and $(\mathcal{P}_3(t)\psi)\in AC^2([0,T];H)$. Furthermore, for the second derivative of $(\mathcal{P}_4F)(t)$, we get
\begin{equation*}
    \partial_{tt}(\mathcal{P}_4F)(t)\qquad\qquad\qquad\qquad\qquad\qquad\qquad\qquad\qquad\qquad\qquad\qquad\qquad\qquad\qquad\qquad\qquad\qquad
\end{equation*}
\begin{equation*}
    =\sum\limits_{n=1}^{\infty}\frac{\lambda_n}{\rho(\lambda_nT^{\alpha})}\Big[E_{\alpha,1}(-\lambda_n T^{\alpha})t^{\alpha-2}E_{\alpha,\alpha-1}(-\lambda_n t^{\alpha})+\lambda_nT^{\alpha-1}E_{\alpha,\alpha}(-\lambda_n T^{\alpha})t^{\alpha-1}E_{\alpha,\alpha}(-\lambda_n t^{\alpha})\Big]
\end{equation*}
\begin{equation*}
    \times F_n(r)\star r^{\alpha-1}E_{\alpha,\alpha}(-\lambda_nr^{\alpha})\Big|_{r=T}e_n
\end{equation*}
\begin{equation*}
    :=\mathcal{I}_{15}F(t)+\mathcal{I}_{16}F(t).\qquad\quad\quad\quad
\end{equation*}
Then, similarly to \eqref{Eq31}, we get
\begin{equation*}
    \|\mathcal{I}_{15}F(t)\|^2=\sum\limits_{n=1}^{\infty}\left|\frac{\lambda_n}{\rho(\lambda_nT^{\alpha})}E_{\alpha,1}(-\lambda_nT^{\alpha})t^{\alpha-2}E_{\alpha,\alpha-1}(-\lambda_nt^{\alpha})\int_0^TF_n(s)(T-s)^{\alpha-1}E_{\alpha,\alpha}(-\lambda_n(T-s)^{\alpha})ds\right|^2
\end{equation*}
\begin{equation*}
    \le c^6c_1^{-2}\sum\limits_{n=1}^{\infty}\frac{\lambda_n^{2}T^{2\alpha}}{(1+\lambda_nT^{\alpha})^2}\left(\frac{t^{\alpha-2}}{1+\lambda_nt^{\alpha}}\right)^2\left|\int_0^T\lambda_nF_n(s)\frac{(T-s)^{\alpha-1}}{1+\lambda_n(T-s)^{\alpha}}ds\right|^2
\end{equation*}
\begin{equation*}
    \le c^6c_1^{-2}t^{2(\alpha-2)}\left|\int_0^T\Big(\sum\limits_{n=1}^{\infty}\Big(AF(s,\cdot),e_n\Big)^2\Big)^{1/2}(T-s)^{\alpha-1}ds\right|^2
\end{equation*}
\begin{equation}\label{Eq47}
    \le c_{21}^2t^{2(\alpha-2)}\|F\|^2_{L^{\infty}(0,T;D(A))}.
\end{equation}
Similarly, we have
\begin{equation}\label{Eq48}
    \|\mathcal{I}_{16}F(t)\|\le c_{22}t^{\alpha-2}\|F\|_{L^{\infty}(0,T;D(A))}.
\end{equation}
Hence, note that $1<\alpha<2$, we have $\partial_{tt}\left(\mathcal{P}_4F(t)\right)\in L^1(0,T;H)$ and $\mathcal{P}_4F(t)\in AC^2([0,T];H)$.

Further, by Lemma 1.3, we have
\begin{equation*}
    \partial_{tt}(\mathcal{P}_5F)(t)\qquad\qquad\qquad\qquad\qquad\qquad\qquad\qquad\qquad\qquad\qquad\qquad\qquad\qquad\qquad\qquad\qquad\qquad
\end{equation*}
\begin{equation*}
    =-\sum\limits_{n=1}^{\infty}\frac{\lambda_n}{\rho(\lambda_nT^{\alpha})}\Big[TE_{\alpha,2}(-\lambda_n T^{\alpha})t^{\alpha-2}E_{\alpha,\alpha-1}(-\lambda_n t^{\alpha})-(\gamma+E_{\alpha,1}(-\lambda_n T^{\alpha})t^{\alpha-1}E_{\alpha,\alpha}(-\lambda_n t^{\alpha})\Big]
\end{equation*}
\begin{equation*}
    \times F_n(r)\star r^{\alpha-2}E_{\alpha,\alpha-1}(-\lambda_nr^{\alpha})\Big|_{r=T}e_n
\end{equation*}
\begin{equation*}
    :=\mathcal{I}_{17}F(t)+\mathcal{I}_{18}F(t).\qquad\quad\quad\quad
\end{equation*}
Then
\begin{equation*}
    \|\mathcal{I}_{17}F(t)\|^2\qquad\qquad\qquad\qquad\qquad\qquad\qquad\qquad\qquad\qquad\qquad\qquad\qquad\qquad\qquad\qquad\qquad\quad\quad
\end{equation*}
\begin{equation*}
    =\sum\limits_{n=1}^{\infty}\left|\frac{\lambda_n}{\rho(\lambda_nT^{\alpha})}TE_{\alpha,2}(-\lambda_nT^{\alpha})t^{\alpha-2}E_{\alpha,\alpha-1}(-\lambda_nt^{\alpha})\int_0^TF_n(s)(T-s)^{\alpha-2}E_{\alpha,\alpha-1}(-\lambda_n(T-s)^{\alpha})ds\right|^2
\end{equation*}
\begin{equation*}
    \le c^6c_1^{-2}\sum\limits_{n=1}^{\infty}\frac{\lambda_n^{2}T^{2\alpha+2}}{(1+\lambda_nT^{\alpha})^2}\left(\frac{t^{\alpha-2}}{1+\lambda_nt^{\alpha}}\right)^2\left|\int_0^T\lambda_nF_n(s)\frac{(T-s)^{\alpha-2}}{1+\lambda_n(T-s)^{\alpha}}ds\right|^2
\end{equation*}
\begin{equation*}
    \le c^6c_1^{-2}t^{2(\alpha-2)}\left|\int_0^T\Big(\sum\limits_{n=1}^{\infty}\Big(AF(s,\cdot),e_n\Big)^2\Big)^{1/2}(T-s)^{\alpha-2}ds\right|^2
\end{equation*}
\begin{equation}\label{Eq49}
    \le c_{23}^2t^{2(\alpha-2)}\|F\|^2_{L^{\infty}(0,T;D(A))}.
\end{equation}
Similarly to \eqref{Eq49}, we obtain
\begin{equation}\label{Eq50}
    \|\mathcal{I}_{18}F(t)\|\le c_{24}t^{\alpha-2}\|F\|_{L^{\infty}(0,T;D(A^{2-\frac{1}{\alpha}}))}.
\end{equation}
Hence, note that $1<\alpha<2$, we have $\partial_{tt}\left(\mathcal{P}_5F(t)\right)\in L^1(0,T;H)$ and $\mathcal{P}_5F(t)\in AC^2([0,T];H)$.

Then, adding the \eqref{Eq44}-\eqref{Eq50} and also \eqref{Eq1}, we obtain
\begin{equation*}
    \partial_t^{\alpha}u(t)=-Au(t)+F(t)\in L^{\infty}(0,T;H).
\end{equation*}
Therefore, $u\in AC^2([0, T]; H)$ and estimate \eqref{Th1.3} is easy to obtain, we argue to complete the proof.

Now, we prove $u\in C^{\theta}([0,T];H)$. Let us consider $0\le t<t+h\le T$. Then 
\begin{equation*}
    u(t+h)-u(t)\qquad\qquad\qquad\qquad\qquad\qquad\qquad\qquad\qquad\qquad\qquad\qquad\qquad
\end{equation*}
\begin{equation*}
    =\Big[(\mathcal{P}_1F)(t+h)-(\mathcal{P}_1F)(t)\Big]+\Big[\mathcal{P}_2(t+h)\varphi-\mathcal{P}_2(t)\varphi\Big]
\end{equation*}
\begin{equation*}
    +\Big[\mathcal{P}_3(t+h)\psi-\mathcal{P}_3(t)\psi\Big]+\Big[(\mathcal{P}_4F)(t+h)-(\mathcal{P}_4F)(t)\Big]
\end{equation*}
\begin{equation*}
    +\Big[(\mathcal{P}_5F)(t+h)-(\mathcal{P}_5F)(t)\Big]:=\mathrm{I}_1+\mathrm{I}_2+\mathrm{I}_3+\mathrm{I}_4+\mathrm{I}_5.\,\,\,
\end{equation*}
By using the differentiation identities
\begin{equation*}
    \partial_tE_{\alpha,1}(-\lambda_nt^{\alpha})=-\lambda_nt^{\alpha-1}E_{\alpha,\alpha}(-\lambda_nt^{\alpha}),\quad \partial_t(t^{\alpha-1}E_{\alpha,\alpha}(-\lambda_nt^{\alpha}))=t^{\alpha-2}E_{\alpha,\alpha-1}(-\lambda_nt^{\alpha}),
\end{equation*}
we have
\begin{equation*}
    \mathrm{I}_1=\sum\limits_{n=1}^{\infty}\Big[\int_0^{t+h}F_n(s)(t+h-s)^{\alpha-1}E_{\alpha,\alpha}(-\lambda_n(t+h-s)^{\alpha})ds\qquad\qquad\qquad\qquad\qquad\qquad
    \end{equation*}
\begin{equation*}
\qquad\qquad\qquad\quad\,-\int_0^{t}F_n(s)(t-s)^{\alpha-1}E_{\alpha,\alpha}(-\lambda_n(t-s)^{\alpha})ds\Big]e_n
\end{equation*}
\begin{equation*}
=\sum\limits_{n=1}^{\infty}\Big[\int_t^{t+h}F_n(s)(t+h-s)^{\alpha-1}E_{\alpha,\alpha}(-\lambda_n(t+h-s)^{\alpha})ds\qquad\qquad\qquad\qquad\qquad\quad
    \end{equation*}
\begin{equation*}
+\int_0^{t}F_n(s)\left((t+h-s)^{\alpha-1}E_{\alpha,\alpha}(-\lambda_n(t+h-s)^{\alpha})-(t-s)^{\alpha-1}E_{\alpha,\alpha}(-\lambda_n(t-s)^{\alpha})\right)ds\Big]e_n
\end{equation*}
\begin{equation*}
=\sum\limits_{n=1}^{\infty}\Big[\int_t^{t+h}F_n(s)(t+h-s)^{\alpha-1}E_{\alpha,\alpha}(-\lambda_n(t+h-s)^{\alpha})ds\Big]e_n
    \end{equation*}
\begin{equation*}
+\sum\limits_{n=1}^{\infty}\Big[\int_0^{t}F_n(s)\tau^{\alpha-1}E_{\alpha,\alpha}(-\lambda_n\tau^{\alpha})\bigg|_{\tau=t-s}^{\tau=t+h-s}ds\Big]e_n
\end{equation*}
or
\begin{equation*}
    \mathrm{I}_1=\sum\limits_{n=1}^{\infty}\Big[\int_t^{t+h}F_n(s)(t+h-s)^{\alpha-1}E_{\alpha,\alpha}(-\lambda_n(t+h-s)^{\alpha})ds\Big]e_n\qquad\qquad\qquad\qquad\qquad\qquad
    \end{equation*}
    \begin{equation}\label{Eq51}
    +\sum\limits_{n=1}^{\infty}\Big[\int_0^{t}\int_{t-s}^{t+h-s}F_n(s)\tau^{\alpha-2}E_{\alpha,\alpha-1}(-\lambda_n\tau^{\alpha})d\tau ds\Big]e_n:=\mathrm{I}_{11}+\mathrm{I}_{12}.
\end{equation}
Now, we establish estimates for $\mathrm{I}_j,j=1,2,$ and show that $\mathrm{I}_j$ tends to $0$ as $h\rightarrow0$. By Lemma 1.1 and \eqref{Eq25}, we have
\begin{equation*}
    \|\mathrm{I}_{11}(t)\|^2_{D(A)}\le c^2\sum\limits_{n=1}^{\infty}\Big|\int_t^{t+h}\lambda_n^{1-\frac{1}{\alpha}}F_n(s)\frac{(\lambda_n(t+h-s)^{\alpha})^{\frac{1}{\alpha}}}{1+\lambda_n(t+h-s)^{\alpha}}(t+h-s)^{\alpha-2}ds\Big|^2\qquad\qquad\qquad\qquad
\end{equation*}
\begin{equation*}
    \le c^2\frac{(\alpha-1)^{\frac{2\alpha-2}{\alpha}}}{\alpha^2}\left(\int_t^{t+h}\Big(\sum\limits_{n=1}^{\infty}\big(\lambda_n^{1-\frac{1}{\alpha}}F_n(s)\big)^2\Big)^{1/2}(t+h-s)^{\alpha-2}ds\right)^2
\end{equation*}
\begin{equation*}
    \le c^2\frac{(\alpha-1)^{\frac{2\alpha-2}{\alpha}}}{\alpha^2}\|F\|^2_{C([0,T];D(A^{1-\frac{1}{\alpha}}))}h^{2(\alpha-1)}.\qquad\qquad\qquad\qquad\qquad\,
\end{equation*}
This leads to 
\begin{equation}\label{Eq52}
    \|\mathrm{I}_{11}(t)\|_{D(A)}\le c_{25}\|F\|_{C([0,T];D(A^{1-\frac{1}{\alpha}}))}h^{\alpha-1},
\end{equation}
where $c_{25}=\frac{(\alpha-1)^{\frac{\alpha-1}{\alpha}}}{\alpha}c$. Secondly, an estimate for the term $\mathrm{I}_{12}$ can be shown by using Lemma 1.1 as follows 
\begin{equation*}
    \|\mathrm{I}_{12}(t)\|^2_{D(A)}\le c^2\sum\limits_{n=1}^{\infty}\Big|\int_0^t\lambda_nF_n(s)\int_{t-s}^{t+h-s}\tau^{\alpha-2}ds\Big|^2.
\end{equation*}
First, for $0<s<t$, we have
\begin{equation*}
    \int_{t-s}^{t+h-s}\tau^{\alpha-2}ds=\frac{1}{\alpha-1}\left[(t+h-s)^{\alpha-1}-(t-s)^{\alpha-1}\right],
\end{equation*}
where we note that the estimate $(t+h-s)^{\alpha-1}-(t-s)^{\alpha-1}\le h^{\alpha-1}$ follows from the properties of the concave function, i.e., $f(t)=t^{\kappa},\,0<\kappa<1: \,f(t_1)-f(t_2)\le f(t_1-t_2)$ if $0<t_2<t_1$. Therefore, we deduce
\begin{equation*}
    \|\mathrm{I}_{12}(t)\|_{D(A)}^2\le \frac{c^2}{(\alpha-1)^2}\sum\limits_{n=1}^{\infty}\Big|\int_0^t\lambda_nF_n(s)ds\Big|^2h^{2\alpha-2}
\end{equation*}
or
\begin{equation}\label{Eq53}
    \|\mathrm{I}_{12}(t)\|_{D(A)}\le c_{26}\|F\|_{C([0,T];D(A))}h^{\alpha-1}.
\end{equation}
Further, 
\begin{equation*}
\mathrm{I}_2=\sum\limits_{n=1}^{\infty}\Big[E_{\alpha,1}(-\lambda_nT^{\alpha})\big(E_{\alpha,1}(-\lambda_n(t+h)^{\alpha})-E_{\alpha,1}(-\lambda_nt^{\alpha})\big)\Big]\frac{\varphi_n}{\rho(\lambda_nT^{\alpha})}e_n\qquad\qquad\qquad\qquad
\end{equation*}
\begin{equation*}
    +\sum\limits_{n=1}^{\infty}\Big[\lambda_nT^{\alpha-1}E_{\alpha,\alpha}(-\lambda_nT^{\alpha})\big((t+h)E_{\alpha,2}(-\lambda_n(t+h)^{\alpha})-tE_{\alpha,2}(-\lambda_nt^{\alpha})\big)\Big]\frac{\varphi_n}{\rho(\lambda_nT^{\alpha})}e_n
\end{equation*}
\begin{equation*}
    :=\mathrm{I}_{21}+\mathrm{I}_{22}.\qquad\qquad\qquad\qquad\qquad\qquad\qquad\qquad\qquad\qquad\qquad\qquad\qquad\qquad\qquad\quad\,\,
\end{equation*}
In order to estimate $\mathrm{I}_{21}$, by Lemma 1.1 and 1.3, we have
\begin{equation*}
    |E_{\alpha,1}(-\lambda_n(t+h)^{\alpha})-E_{\alpha,1}(-\lambda_nt^{\alpha})|=\left|\int_{t+h}^t\lambda_n\tau^{\alpha-1}E_{\alpha,\alpha}(-\lambda_n\tau^{\alpha})d\tau\right|\qquad\qquad\qquad\qquad
\end{equation*}
\begin{equation}\label{Eq54}
    \le\int_t^{t+h}\lambda_n\tau^{\alpha-1}\frac{c}{1+\lambda_n\tau^{\alpha}}d\tau\le c\int_t^{t+h}\tau^{-1}d\tau.
\end{equation}
Then for $\delta\le t\le T$, we have
\begin{equation*}
    \|\mathrm{I}_{21}(t)\|^2_{D(A)}\le c^2c_1^{-2}\sum\limits_{n=1}^{\infty}
    \left|\frac{\lambda_n^2T^{\alpha}}{1+\lambda_nT^{\alpha}}\varphi_n\int_{t}^{t+h}\tau^{-1}d\tau\right|^2\le c^2c_1^{-2}\left(\ln\Big(1+\frac{h}{t}\Big)\right)^2\sum\limits_{n=1}^{\infty}\lambda_n^2\varphi_n^2
\end{equation*}
\begin{equation}\label{Eq55}
    \le c^2c_1^{-2}\frac{h^2}{\delta^2}\|\varphi\|_{D(A)}^2\le c_{27}^2\frac{h^{2(\alpha-1)}}{\delta^2}\|\varphi\|_{D(A)}^2.
\end{equation}
Here we use $\ln(1+\eta)\le \eta$ for $\eta>0$. Now we will estimate $\mathrm{I}_{22}$.  Herein, by Lemma 1.1 and 1.3, we have
\begin{equation*}
\|\mathrm{I}_{22}(t)\|^2_{D(A)}\le c^2c_1^{-2}\sum\limits_{n=1}^{\infty}\left|\frac{\lambda_n^2T^{2\alpha-1}}{1+\lambda_nT^{\alpha}}\varphi_n\int_{t}^{t+h}E_{\alpha,1}(-\lambda_ns^{\alpha})ds\right|^2
\end{equation*}
\begin{equation*}
    \le c^4c_1^{-2}T^{2(\alpha-1)}\sum\limits_{n=1}^{\infty}\varphi_n^2\left|\int_t^{t+h}s^{-\alpha}ds\right|^2\le \frac{c^4c_1^{-2}}{(\alpha-1)^2}T^{2(\alpha-1)}\left|\frac{(t+h)^{\alpha-1}-t^{\alpha-1}}{(t+h)^{\alpha-1}t^{\alpha-1}}\right|^2\|\varphi\|^2.
\end{equation*}
We note $(t+h)^{\alpha-1}-t^{\alpha-1}\le h^{\alpha-1}$, and for $\delta\le t\le T$ the estimate $(t+h)^{\alpha-1}t^{\alpha-1}\ge\delta^{2\alpha-2}$ holds. Hence, we get
\begin{equation}\label{Eq56}
    \|\mathrm{I}_{22}(t)\|_{D(A)}\le c_{28}\frac{h^{\alpha-1}}{\delta^{2(\alpha-1)}}\|\varphi\|.
\end{equation}

Similar to \eqref{Eq54}-\eqref{Eq56}, the following estimate is hold for $\mathrm{I}_3$:
\begin{equation}\label{Eq57}
    \|\mathrm{I}_3(t,\cdot)\|_{D(A)}\le c_{29}\frac{h^{\alpha-1}}{\delta}\|\psi\|+c_{30}\frac{h^{\alpha-1}}{\delta^{2(\alpha-1)}}\|\psi\|_{D(A)}.
\end{equation}
We will now estimate $\mathrm{I}_4$. First, we have
\begin{equation*}
    \mathrm{I}_4=-\sum\limits_{n=1}^{\infty}\Big[E_{\alpha,1}(-\lambda_nT^{\alpha})\big(E_{\alpha,1}(-\lambda_n(t+h)^{\alpha})-E_{\alpha,1}(-\lambda_nt^{\alpha})\big)\qquad\qquad\qquad\qquad\qquad\qquad
\end{equation*}
\begin{equation*}
    +\big((t+h)E_{\alpha,2}(-\lambda_n(t+h)^{\alpha})-tE_{\alpha,2}(-\lambda_nt^{\alpha})\big)\Big]\frac{F_n(r)\star r^{\alpha-1}E_{\alpha,\alpha}(-\lambda_nr^{\alpha})}{\rho(\lambda_nT^{\alpha})}e_n
\end{equation*}
\begin{equation*}
    :=\mathrm{I}_{41}+\mathrm{I}_{42}.
\end{equation*}
Since the main part of $\mathrm{I}_4$ is the same as $\mathrm{I}_2$, therefore, we write it directly without calculating. First, for $\mathrm{I}_{41}$, we have
\begin{equation*}
    \|\mathrm{I}_{41}(t)\|_{D(A)}^2\le c_{27}^2\frac{h^{2(\alpha-1)}}{\delta^2}\sum\limits_{n=1}^{\infty}\left|\frac{\lambda_nT^{\alpha}}{1+\lambda_nT^{\alpha}}\int_0^T\lambda_n^{1-\frac{1}{\alpha}}F_n(s)\frac{(\lambda_n(T-s)^{\alpha})^{\frac{1}{\alpha}}}{1+\lambda_n(T-s)^{\alpha}}(T-s)^{\alpha-2}ds\right|^2
\end{equation*}
\begin{equation*}
    \le c_{27}^2\frac{(\alpha-1)^{\frac{2\alpha-2}{\alpha}}}{\alpha^2}\frac{h^{2(\alpha-1)}}{\delta^2}\left|\int_0^T\Big(\sum\limits_{n=1}^{\infty}\Big(\lambda_n^{1-\frac{1}{\alpha}}F_n(s)\Big)^2\Big)^{1/2}(T-s)^{\alpha-2}ds\right|^2
\end{equation*}
\begin{equation}\label{Eq58}
    \le c_{31}^2\frac{h^{2(\alpha-1)}}{\delta^2}\|F\|_{C([0,T];D(A^{1-\frac{1}{\alpha}}))}^2.
\end{equation}
Similarly, we can show
\begin{equation}\label{Eq59}
    \|\mathrm{I}_{42}(t)\|_{D(A)}\le c_{32} \frac{h^{\alpha-1}}{\delta^{2(\alpha-1)}}\|F\|_{C([0,T];L^2(\Omega))},
\end{equation}
and we can prove
\begin{equation}\label{Eq60}
    \|\mathrm{I}_{5}(t)\|_{D(A)}\le c_{33}\frac{h^{\alpha-1}}{\delta}\|F\|_{C([0,T];D(A))}+c_{34}\frac{h^{\alpha-1}}{\delta^{2(\alpha-1)}}\|F\|_{C([0,T];D(A))}.
\end{equation}

Now we complete the proof of Theorem 2.1(4). By Lemma 1.2 and 1.3, we have
\begin{equation*}
    \partial_t^{\alpha}u(t)=\sum\limits_{n=1}^{\infty}F_n(t)e_n  -\sum\limits_{n=1}^{\infty}\lambda_nF_n(t)\star t^{\alpha-1}E_{\alpha,\alpha}(-\lambda_n t^{\alpha})e_n\qquad\qquad\qquad\qquad\qquad\qquad\qquad
\end{equation*}
\begin{equation*}
    -\sum\limits_{n=1}^{\infty}\lambda_n\Big[E_{\alpha,1}(-\lambda_n T^{\alpha})E_{\alpha,1}(-\lambda_n t^{\alpha})+\lambda_nT^{\alpha-1}E_{\alpha,\alpha}(-\lambda_n T^{\alpha})tE_{\alpha,2}(-\lambda_n t^{\alpha})\Big]\frac{\varphi_n}{\rho(\lambda_nT^{\alpha})}e_n
\end{equation*}
\begin{equation*}
    \quad-\sum\limits_{n=1}^{\infty}\lambda_n\Big[-TE_{\alpha,2}(-\lambda_n T^{\alpha})E_{\alpha,1}(-\lambda_n t^{\alpha})+(\gamma+E_{\alpha,1}(-\lambda_nT^{\alpha}))tE_{\alpha,2}(-\lambda_n t^{\alpha})\Big]\frac{\psi_n}{\rho(\lambda_nT^{\alpha})}e_n
\end{equation*}
\begin{equation*}
    +\sum\limits_{n=1}^{\infty}\lambda_n\Big[E_{\alpha,1}(-\lambda_n T^{\alpha})E_{\alpha,1}(-\lambda_n t^{\alpha})+\lambda_nT^{\alpha-1}E_{\alpha,\alpha}(-\lambda_n T^{\alpha})tE_{\alpha,2}(-\lambda_n t^{\alpha})\Big]\qquad\qquad\qquad
\end{equation*}
\begin{equation*}
    \times\frac{F_n(r)\star r^{\alpha-1}E_{\alpha,\alpha}(-\lambda_nr^{\alpha})\Big|_{r=T}}{\rho(\lambda_nT^{\alpha})}e_n
\end{equation*}
\begin{equation*}
    -\sum\limits_{n=1}^{\infty}\lambda_n\Big[TE_{\alpha,2}(-\lambda_n T^{\alpha})E_{\alpha,1}(-\lambda_n t^{\alpha})-(\gamma+E_{\alpha,1}(-\lambda_nT^{\alpha}))tE_{\alpha,2}(-\lambda_n t^{\alpha})\Big]\qquad\qquad\qquad
\end{equation*}
\begin{equation*}
    \times\frac{F_n(r)\star r^{\alpha-2}E_{\alpha,\alpha-1}(-\lambda_nr^{\alpha})\Big|_{r=T}}{\rho(\lambda_nT^{\alpha})}e_n.
\end{equation*}
Then
\begin{equation*}
    \partial_t^{\alpha}u(t+h)-\partial_t^{\alpha}u(t)\qquad\qquad\qquad\qquad\qquad\qquad\qquad\qquad\qquad\qquad\qquad\qquad\qquad\qquad\qquad\qquad
\end{equation*}
\begin{equation*}
    =F(t+h)-F(t)-\Big[(\mathcal{P}_1AF)(t+h)-(\mathcal{P}_1AF)(t)\Big]\qquad\qquad\qquad\quad
\end{equation*}
\begin{equation*}
    -\Big[(\mathcal{P}_2(t+h)A\varphi)-(\mathcal{P}_2(t)A\varphi)\Big] -\Big[(\mathcal{P}_3(t+h)A\psi)-(\mathcal{P}_3(t)A\psi)\Big]\qquad
\end{equation*}
\begin{equation*}
   +\Big[(\mathcal{P}_4AF)(t+h)-(\mathcal{P}_4AF)(t)\Big]
    -\Big[(\mathcal{P}_5AF)(t+h)-(\mathcal{P}_5AF)(t)\Big]\,\,\,\,
\end{equation*}
\begin{equation}\label{Eq61}
    :=\mathrm{I}_6+\mathrm{I}_7+\mathrm{I}_8+\mathrm{I}_9+\mathrm{I}_{10}+\mathrm{I}_{11}.\qquad\qquad\qquad\qquad\qquad\qquad\qquad\qquad\qquad\quad
\end{equation}
Now, we establish estimates for $\mathrm{I}_j\,(j=6,...,11),$ and show that $\mathrm{I}_j$ tends to $0$ as $h\rightarrow0$. The next estimate can be obtained easily
\begin{equation*}
    \|\mathrm{I}_6(t)\|=\frac{\|F(t+h)-F(t)\|}{h^{\alpha-1}}h^{\alpha-1}\le[F]_{C^{0,\alpha-1}([0,T];H)}h^{\alpha-1}.
\end{equation*}
By \eqref{Eq51}, the second term of \eqref{Eq61} may write as follows
\begin{equation*}
    \mathrm{I}_7=-\sum\limits_{n=1}^{\infty}\lambda_n\Big[\int_t^{t+h}F_n(s)(t+h-s)^{\alpha-1}E_{\alpha,\alpha}(-\lambda_n(t+h-s)^{\alpha})ds\Big]e_n\qquad\qquad\qquad\qquad\qquad
    \end{equation*}
    \begin{equation}\label{Eq62}
    \qquad+\sum\limits_{n=1}^{\infty}\lambda_n\Big[\int_0^{t}\int_{t-s}^{t+h-s}F_n(s)\tau^{\alpha-2}E_{\alpha,\alpha-1}(-\lambda_n\tau^{\alpha})d\tau ds\Big]e_n.
\end{equation}
We can see that the last equality differs from the \eqref{Eq51} only by the factor $-\lambda_n$ value. Therefore, we obtain
\begin{equation}\label{Eq63}
    \|\mathrm{I}_7(t)\|\le c_{35}\left(\|F\|_{C([0,T];D(A^{1-\frac{1}{\alpha}}))}+\|F\|_{C([0,T];D(A))}\right)h^{\alpha-1}.
\end{equation}
Similarly to \eqref{Eq55} and \eqref{Eq56}, we have
\begin{equation}\label{Eq64}
    \|\mathrm{I}_8(t)\|\le c_{36}\left(\frac{1}{\delta}\|\varphi\|_{D(A)}+\frac{1}{\delta^{\alpha-1}}\|\varphi\|\right)h^{\alpha-1}.
\end{equation}
Further, by \eqref{Eq57}, we obtain
\begin{equation}\label{Eq65}
    \|\mathrm{I}_9(t)\|\le c_{37}\left(\frac{1}{\delta}\|\varphi\|_{D(A)}+\frac{1}{\delta^{\alpha-1}}\|\varphi\|\right)h^{\alpha-1}.
\end{equation}
The same
\begin{equation}\label{Eq66}
    \|\mathrm{I}_{10}(t)\|\le c_{38}\left(\frac{1}{\delta}\|F\|_{C([0,T];D(A^{1-\frac{1}{\alpha}})}+\frac{1}{\delta^{\alpha-1}}\|F\|_{C([0,T];H)}\right)h^{\alpha-1}
\end{equation}
and
\begin{equation}\label{Eq67}
    \|\mathrm{I}_{11}(t)\|\le c_{39}\left(\frac{1}{\delta}\|F\|_{C([0,T];D(A)}+\frac{1}{\delta^{\alpha-1}}\|F\|_{C([0,T];D(A))}\right)h^{\alpha-1}.
\end{equation}
Thus the proof of Theorem 2.1 is complete.

{\bf Corollary 1.} \textit{Let $\gamma\not=0$ and  $T\not\in\Lambda$. Let $\varphi,\psi\in D(A)$ and $F=0$. Then, there exists a constant $C_5>0$ such that $u\in C^{\infty}((0,T];H)$ and
\begin{equation}\label{Eq68}
    \|\partial_t^{m}u(t)\|\le \frac{C_6}{t^m}\Big(\|\varphi\|+\|\psi\|+t(\|\varphi\|_{D(A)}+\|\psi\|_{D(A)})\Big),\quad t>0,\,\,m\in\mathbb{N}.
\end{equation}
}

The proof of the result is a direct consequence of the following relations 
\begin{equation*}
    \frac{d^m}{dt^m}E_{\alpha,1}(-\lambda t^{\alpha})=-\lambda t^{\alpha-m}E_{\alpha,\alpha-m+1}(-\lambda t^{\alpha}),\quad \lambda>0,\,\,t>0
\end{equation*}
and
\begin{equation*}
    \frac{d}{dt} tE_{\alpha,2}(-\lambda t^{\alpha})=E_{\alpha,1}(-\lambda t^{\alpha}),\quad t\ge0.
\end{equation*}

\section{Inverse problem 1}

\ \ \ \ We set
\begin{equation*}
\mathrm{P}_1=\{\xi>0:\,\bigtriangleup_n(\xi)=0\},
\end{equation*}
where
\begin{equation*}
    \bigtriangleup_n(t):=t^{\alpha}E_{\alpha,\alpha+1}(-\lambda_nt^{\alpha})+T^{\alpha}E_{\alpha,\alpha+1}(-\lambda_nT^{\alpha})\mathcal{P}_{n}^{(2)}(t)+T^{\alpha-1}E_{\alpha,\alpha}(-\lambda_nT^{\alpha})\mathcal{P}_{n}^{(3)}(t),
\end{equation*}
\begin{equation*}
    \mathcal{P}_{n}^{(2)}(t)=\Big[E_{\alpha,1}(-\lambda_n T^{\alpha})E_{\alpha,1}(-\lambda_n t^{\alpha})+\lambda_nT^{\alpha-1}E_{\alpha,\alpha}(-\lambda_n T^{\alpha})tE_{\alpha,2}(-\lambda_n t^{\alpha})\Big]\frac{1}{\rho(\lambda_nT^{\alpha})},
\end{equation*}
and
\begin{equation*}
    \mathcal{P}_{n}^{(3)}(t)=\Big[-TE_{\alpha,2}(-\lambda_n T^{\alpha})E_{\alpha,1}(-\lambda_n t^{\alpha})+(\gamma+E_{\alpha,1}(-\lambda_nT^{\alpha}))tE_{\alpha,2}(-\lambda_n t^{\alpha})\Big]\frac{1}{\rho(\lambda_nT^{\alpha})},
\end{equation*}
here $\rho(\eta)$ is given by \eqref{Eq13}.

{\bf Theorem 3.1.} \textit{Let $\gamma\not=0$, $T\not\in\Lambda,\xi\not\in\mathrm{P}_1$. If $\varphi,\psi\in D(A^{2-\frac{1}{\alpha}})$ and $h\in D(A)$, then the inverse problem \eqref{Eq1}-\eqref{Eq3} has a unique solution $\{u(t),f\}$ and this solution has the following form
\begin{equation}\label{Eq69}
u(t)=\mathcal{P}_2(t)\varphi+\mathcal{P}_3(t)\psi+\sum\limits_{n=1}^{\infty}\bigtriangleup_n(t)\cdot f_ne_n,
    \end{equation}
\begin{equation}\label{Eq70}
    f=\sum\limits_{n=1}^{\infty}\left(h_n-\mathcal{P}_{n}^{(2)}(\xi)\cdot\varphi_n-\mathcal{P}_{n}^{(3)}(\xi)\cdot\psi_n\right)\frac{e_n}{\bigtriangleup_n(\xi)},
\end{equation}
where $f_n=(f,e_n),\,h_n=(h,e_n)$.
}

{\bf Proof.}  If $f$ is known, then the unique solution of problem \eqref{Eq1} has the form \eqref{Eq23}, and since $f$ does not depend on $t$, then, thanks to formula
\begin{equation}\label{Eq71}
    \int_0^ts^{\beta-1}E_{\alpha,\beta}(-\lambda s^{\alpha})ds=t^{\beta}E_{\alpha,\beta+1}(-\lambda t^{\alpha}),\quad \beta>0
\end{equation}
(see \cite{Podlubny}, p.24), it is easy to verify that the formal solution of problem \eqref{Eq1}-\eqref{Eq2} has the form \eqref{Eq69}.

By virtue of additional condition \eqref{Eq3} and completeness of the system $e_n$, and furthermore from \eqref{Eq23}, we obtain:
\begin{equation*}
    \bigg\{\int_0^{\xi}s^{\alpha-1}E_{\alpha,\alpha}(-\lambda_ns^{\alpha})ds-\int_0^{T}s^{\alpha-1}E_{\alpha,\alpha}(-\lambda_ns^{\alpha})ds\cdot\mathcal{P}_n^{(2)}(\xi)\qquad\qquad\qquad\qquad\qquad\qquad
\end{equation*}
\begin{equation*}
    \qquad\qquad\qquad-\int_0^{T}s^{\alpha-2}E_{\alpha,\alpha-1}(-\lambda_ns^{\alpha})ds\cdot\mathcal{P}_n^{(3)}(\xi)\bigg\}f_n=h_n-\varphi_n\cdot\mathcal{P}_n^{(2)}(\xi)-\psi_n\cdot\mathcal{P}_n^{(3)}(\xi).
\end{equation*}
Using the formula \eqref{Eq71}, and after simple calculations, we arrive at \eqref{Eq70}. We emphasise that the main problem is determining whether $\bigtriangleup_n(\xi)$ has zeros.  For this, first, we give the asymptotic behaviour of the $E_{\alpha,\alpha+1}(-\eta)$ for $1<\alpha<2$ (see, \cite{Podlubny}, p.34), by
\begin{equation}\label{Eq72}
    E_{\alpha,\alpha+1}(-\eta)=\frac{1}{\eta}+O\left(\frac{1}{\eta^2}\right),\quad\mbox{as}\quad \eta\rightarrow\infty.
\end{equation}
Using the asymptotic expansions \eqref{Eq15.0} and \eqref{Eq72}, we have
\begin{equation}\label{Eq72.0}
    \mathcal{P}_n^{(2)}(\xi)=\left(\frac{1}{\alpha-1}\frac{\xi}{T}-\frac{1}{\alpha}\right)\frac{1}{\gamma\xi^{\alpha}(\alpha-1)}\frac{1}{\lambda_n}+O\left(\frac{1}{\lambda_n^2}\right),
\end{equation}
and
\begin{equation}\label{Eq72.1}
    \mathcal{P}_n^{(3)}(\xi)=\frac{1}{1-\alpha}\frac{T^{\alpha}}{\xi^{\alpha-1}}+O\left(\frac{1}{\lambda_n}\right).
\end{equation}
By applying \eqref{Eq72.0} and \eqref{Eq72.1}, we obtain the following result for $\bigtriangleup_n(\xi)$:
\begin{equation}\label{Eq72.2}
    \bigtriangleup_n(\xi)=\frac{1}{\lambda_n}+O\left(\frac{1}{\lambda_n^2}\right).
\end{equation}
On the other hand, we have
\begin{equation}\label{Eq72.3}
     \bigtriangleup_n(0)=\frac{\gamma^{-1}}{\lambda_n}+O\left(\frac{1}{\lambda_n^2}\right).
\end{equation}
Therefore, the same reason as Lemma 2.1, the set $\mathrm{P}_1$ (for $\gamma\not=0$) is a non-empty and finite set. Further, our next aim will get estimates of $\mathcal{P}_n^{(2)}(\xi),\mathcal{P}_n^{(2)}(\xi)$ and $\bigtriangleup_n(\xi)$ from above and below respectively. Therefore, by 
\eqref{Eq72.0} and \eqref{Eq72.1}, we can choose sufficiently a large constant $N_0>0$ such that
\begin{equation}\label{Eq73}
    \frac{c_{40}}{\lambda_n}\le\left|\mathcal{P}_n^{(2)}(\xi)\right|\le\frac{c_{41}}{\lambda_n},\quad n\ge N_0,
\end{equation}
and
\begin{equation}\label{Eq74}
    \left|\mathcal{P}_n^{(3)}(\xi)\right|\le c_{42},\quad n\ge N_0,
\end{equation}
here  and henceforth  $c_k$ are independent of $n$, but dependent on $T,\xi,N_0$ and $\alpha,\gamma$. By employing the inequalities \eqref{Eq73} and \eqref{Eq74}, we obtain the following result for $\bigtriangleup_n(\xi)$:
\begin{equation}\label{Eq75}
    |\bigtriangleup_n(\xi)|\ge\frac{1}{\lambda_n},\quad \mbox{as}\quad n\rightarrow\infty.
\end{equation}

Let us show the convergence of series \eqref{Eq70}. Indeed, by Parseval's identity and  \eqref{Eq73}-\eqref{Eq75}, we have
\begin{equation*}
    \|f\|^2=\sum\limits_{n=1}^{\infty}\left|\Big(h_n-\mathcal{P}_{n}^{(2)}(\xi)\cdot\varphi_n-\mathcal{P}_{n}^{(3)}(\xi)\cdot\psi_n\Big)\frac{1}{\bigtriangleup_n(\xi)}\right|^2\qquad\qquad\qquad
\end{equation*}
\begin{equation*}
    \le 2\sum\limits_{n=1}^{\infty}\lambda_n^2h_n^2+2c_{41}^2\sum\limits_{n=1}^{\infty}\varphi_n^2+2c_{42}^2\sum\limits_{n=1}^{\infty}\lambda_n^2\psi_n^2\qquad\quad
\end{equation*}
\begin{equation}\label{Eq76}
    \le c_{43}^2\Big(\|h\|_{D(A)}^2+\|\varphi\|^2+\|\psi\|_{D(A)}^2\Big).
\end{equation}
Therefore, $f\in H$. Moreover,  if we make the same calculations as for \eqref{Eq26}, \eqref{Eq27}, and using \eqref{Eq75}, we obtain 
\begin{equation*}
    \|u(t)\|_{D(A)}^2\le\|\mathcal{P}_2(t)A\varphi\|^2+\|\mathcal{P}_3(t)A\psi\|^2+\sum\limits_{n=1}^{\infty}|\bigtriangleup_n(t)\cdot \lambda_nf_n|^2\qquad\qquad\qquad\qquad\qquad
\end{equation*}
\begin{equation}\label{Eq77}
    \le c_{44}^2\|\varphi\|^2_{D(A^{2-\frac{1}{\alpha}})}+c_{45}^2\|\psi\|^2_{D(A^{2-\frac{1}{\alpha}})}+\|f\|^2,\quad t\in[0,T].
\end{equation}
The other conditions of Definition 1.1 will be shown as condition 2) of Theorem 2.1.

{\it Uniqueness.} Suppose we have two solutions: $\{u_1(t),f_1\}$
 and $\{u_2(t),f_2\}.$ It is required to prove $u(t)=u_1(t)-u_2(t)\equiv0$
 and $f=f_1-f_2\equiv0$. Since the problem is linear, hence, to determine $u(t)$
 and $f$ we have the problem:
 \begin{equation}\label{Eq77.0}
    \begin{cases}
        \partial_{t}^{\alpha}u(t)+Au(t)=f,&0<t<T,\\
         u(T)=\partial_tu(T)=0,\\
         u(\xi)=0,
    \end{cases}
\end{equation}
where $\xi\in(0,T)$ fixed number.
Let $u(t)$ be a solution to this problem and $u_n(t)=(u(t),e_n)$. Then, as \eqref{un2}, by virtue of problem \eqref{Eq77.0} and the self-adjointness of operator $A$, we obtain
\begin{equation*}
    \begin{cases}
\partial_{t}^{\alpha}u_n(t)+\lambda_nu_n(t)=f_n,&0<t<T,\\
         \gamma u_n(0)+ u_n(T)=0,\quad \partial_tu_n(T)=0,\\
         u_n(\xi)=0.
    \end{cases}
\end{equation*}
Suppose that $f_n$
 is known and $T\not\in\Lambda$. Then, use the final time conditions to obtain
 \begin{equation*}
     u_n(t)=f_n\bigtriangleup_n(t).
 \end{equation*}
Now apply $u_n(\xi)=0$ to obtain
\begin{equation}\label{Eq77.1}
    f_n\bigtriangleup_n(\xi)=0.
\end{equation}
By $\xi\not\in\mathrm{P}_1$, we have $f_n=0$. Hence, from the completeness of the system of eigenfunctions $e_n$, we finally obtain $f=0$
 and $u(t)\equiv0$, as required. The uniqueness is proved.

\section{Inverse problem 2}

\ \ \ \ In this section, we consider the inverse problem \eqref{Eq1}-\eqref{Eq2}, \eqref{Eq4} for the case $F(t)=fp(t)$, and the following result of stability in the recovery of the time-dependent factor appearing in the source term problem \eqref{Eq1}-\eqref{Eq2}. Therefore, first, we assume the following assumptions:

(C1) $\Phi:\,\big\{\Phi[e_k]\big\}\in l^2(\mathbf{N})$, where $\mathbf{N}$ is a natural number set;

(C2) $\Phi[f]\not=0$;

(C3) $h(t)\in AC[0,T]$;

(C4) $f\in D(A^{\frac{1}{\alpha}})$.

\textbf{Remark 4.1.}  (C1) is not empty. Indeed, let $X=L^2(0,1)$ and $Au=-u_{xx},\,x\in(0,1)$ with homogeneous Dirichlet boundary condition. Then the operator has the eigensystem $\{(\pi k)^2,\sqrt{2}\sin(\pi kx)\}_{k\in\mathbf{N}}$. Let $\Phi[u(t,\cdot)]:=\int_0^1u(t,x)dx$. Since
\begin{equation*}
    \Phi[e_k]=\sqrt{2}\int_0^1\sin(\pi k x)dx=\frac{\sqrt{2}(1+(-1)^{k+1})}{\pi k}=
    \begin{cases}
        \frac{2\sqrt{2}}{\pi k},&k=2m-1\,(m\in\mathbf{N});\\
        0,&k=2m\,(m\in\mathbf{N}),
    \end{cases}
\end{equation*}
we have
\begin{equation*}
    \sum\limits_{k=1}^{\infty}|\Phi[e_k]|^2<\infty.
\end{equation*}

Let $\mathcal{A}u=-u_{xx}$ with homogeneous Neumann boundary condition. In this case, we take the operator $A$ as 
\begin{equation*}
    Au=\mathcal{A}u+u,\quad u\in D(A):=\{u\in H^2(0,1):\,u_x(0,t)=u_x(1,t)=0\}.
\end{equation*}
Then the operator $A$ has the eigensystem $\{1,1\}\cup\{(\pi n)^2+1,\sqrt{2}\cos(\pi nx)\}_{k\in\mathbf{N}}$. Let $\Phi[u(t,\cdot)]:=\int_0^1x(1-x)u(t,x)dx$. Based on the previous calculations, the following results have been obtained
\begin{equation*}
    \Phi[e_k]=
    \begin{cases}
        \frac{1}{6},&k=0;\\
        0,&k=2m-1\,(m\in\mathbf{N});\\
        -\frac{2\sqrt{2}}{(\pi k)^2},&k=2m\,(m\in\mathbf{N}).
    \end{cases}
\end{equation*}
Therefore, we get
\begin{equation*}
    \sum\limits_{k=1}^{\infty}|\Phi[e_k]|^2<\infty.
\end{equation*}

\textbf{Remark 4.2.}  In (C3) implies that the Caputo left derivative $\partial_t^{\alpha}h(t)$ exist almost everywhere on $(0,T]$ (see \cite{Kilbas}, Theorem 2.1, p. 92).

The main result of this section is presented below.

{\bf Theorem 4.1.} \textit{Let $\gamma\not=0$ and $T\not\in\Lambda$. Let (C1)-(C4) hold and let $u$ satisfy \eqref{Eq1}-\eqref{Eq2}, \eqref{Eq4} for $p\in AC[0,T]$. Then there exist  constants $C,C'>0$ such that
\begin{equation}\label{EqTh4.1}
    C\|\partial_t^{\alpha}h\|_{L^{\infty}(0,T)}\le\|p\|_{L^{\infty}(0,T)}\le C'\|\partial_t^{\alpha}h\|_{L^{\infty}(0,T)}.
\end{equation}
}

{\bf Proof.} By $p\in AC[0,T]$ and $f\in D(A^{\frac{1}{\alpha}})$, we apply Theorem 2.1 to obtain
\begin{equation*}
    u(t)=\sum\limits_{n=1}^{\infty}f_n\Big(p(t)\star t^{\alpha-1}E_{\alpha,\alpha}(-\lambda_n t^{\alpha})\Big)e_n-\sum\limits_{n=1}^{\infty}f_n\mathcal{P}_n^{(2)}(t)\cdot\Big(p(r)\star r^{\alpha-1}E_{\alpha,\alpha}(-\lambda_nr^{\alpha})\Big|_{r=T}\Big)e_n
\end{equation*}
\begin{equation*}
    -\sum\limits_{n=1}^{\infty}f_n\mathcal{P}_n^{(3)}(t)\cdot\Big(p(r)\star r^{\alpha-2}E_{\alpha,\alpha-1}(-\lambda_nr^{\alpha})\Big|_{r=T}\Big)e_n
\end{equation*}
in $AC^2([0,T];H)\cap C([0,T];D(A))$ and
\begin{equation*}
    \partial_t^{\alpha}u(t)=fp(t)-\sum\limits_{n=1}^{\infty}\lambda_nf_n \Big(p(t)\star t^{\alpha-1}E_{\alpha,\alpha}(-\lambda_nt^{\alpha})\Big)e_n\qquad\qquad\qquad\qquad\qquad
\end{equation*}
\begin{equation*}
    +\sum\limits_{n=1}^{\infty}\lambda_n\mathcal{P}_n^{(2)}(t)\cdot f_n \Big(p(r)\star r^{\alpha-1}E_{\alpha,\alpha}(-\lambda_nr^{\alpha})\Big)\Big|_{r=T}e_n
\end{equation*}
\begin{equation*}
    \quad+\sum\limits_{n=1}^{\infty}\lambda_n\mathcal{P}_n^{(3)}(t)\cdot f_n \Big(p(r)\star r^{\alpha-2}E_{\alpha,\alpha-1}(-\lambda_nr^{\alpha})\Big)\Big|_{r=T}e_n,
\end{equation*}
in $L^{\infty}(0,T;H)$, where $\mathcal{P}_n^{(2)}(t),\mathcal{P}_n^{(3)}(t)$ are given in Section 3. Moreover, by \eqref{Eq4}, we have
\begin{equation*}
    \partial_t^{\alpha}h(t)=\Phi[f]p(t)-\sum\limits_{n=1}^{\infty}\lambda_nf_n \Big(p(t)\star t^{\alpha-1}E_{\alpha,\alpha}(-\lambda_nt^{\alpha})\Big)\Phi[e_n]\qquad\qquad\qquad\qquad\qquad
\end{equation*}
\begin{equation*}
    +\sum\limits_{n=1}^{\infty}\lambda_n\mathcal{P}_n^{(2)}(t)\cdot f_n \Big(p(r)\star r^{\alpha-1}E_{\alpha,\alpha}(-\lambda_nr^{\alpha})\Big)\Big|_{r=T}\Phi[e_n]
\end{equation*}
\begin{equation}\label{Eq78.0}
    \quad+\sum\limits_{n=1}^{\infty}\lambda_n\mathcal{P}_n^{(3)}(t)\cdot f_n \Big(p(r)\star r^{\alpha-2}E_{\alpha,\alpha-1}(-\lambda_nr^{\alpha})\Big)\Big|_{r=T}\Phi[e_n],
\end{equation}

Now, we will estimate each term of \eqref{Eq78.0}. Therefore, by (C1), \eqref{Eq25}, \eqref{Eq73} and \eqref{Eq74}, we have 
\begin{equation*}
   1) \quad \sum\limits_{n=1}^{\infty}\left|\lambda_nf_n\Big(\int_0^tp(s)(t-s)^{\alpha-1}E_{\alpha,\alpha}(-\lambda_n(t-s)^{\alpha})ds\Big)\Phi[e_n]\right|\qquad\qquad\qquad\qquad\qquad\qquad\quad
\end{equation*}
\begin{equation*}
     \le c\sum\limits_{n=1}^{\infty}\Big|\lambda_n^{\frac{1}{\alpha}}f_n\int_0^tp(s)\frac{(\lambda_n(t-s)^{\alpha})^{\frac{\alpha-1}{\alpha}}}{1+\lambda_n(t-s)^{\alpha}}ds\Phi[e_n]\Big|\qquad\qquad\qquad\quad\,
\end{equation*}
\begin{equation*}
    \le c\frac{(\alpha-1)^{\frac{\alpha-1}{\alpha}}}{\alpha}\sum\limits_{n=1}^{\infty}\Big|\lambda_n^{\frac{1}{\alpha}}f_n\int_0^tp(s)ds\Phi[e_n]\Big|\qquad\qquad\qquad\qquad\quad\,\,\,
\end{equation*}
\begin{equation*}
    \le c\frac{(\alpha-1)^{\frac{\alpha-1}{\alpha}}}{\alpha}t\|p\|_{L^{\infty}(0,T)}\left(\sum\limits_{n=1}^{\infty}\Big(\lambda_n^{\frac{1}{\alpha}}f_n\Big)^2\right)^{\frac{1}{2}}\left(\sum\limits_{n=1}^{\infty}\Big(\Phi[e_n]\Big)^2\right)^{\frac{1}{2}}
\end{equation*}
\begin{equation}\label{Eq78}
    \le c_{46}\|p\|_{L^{\infty}(0,T)}\|f\|_{D(A^{\frac{1}{\alpha}})},\qquad\qquad\qquad\qquad\qquad\qquad\qquad\quad\quad
\end{equation}
and
\begin{equation}\label{Eq78.1}
   2) \quad \sum\limits_{n=1}^{\infty}\bigg|\lambda_n\mathcal{P}_n^{(2)}(t)f_n\int_0^Tp(s)(T-s)^{\alpha-1}E_{\alpha,\alpha}(-\lambda_n(T-s)^{\alpha})ds\bigg|\le c_{47}
\|p\|_{L^{\infty}(0,T)}\|f\|,\qquad\quad
\end{equation}
and for the last term, we obtain
\begin{equation}\label{Eq78.1}
   3) \quad \sum\limits_{n=1}^{\infty}\bigg|\lambda_n\mathcal{P}_n^{(3)}(t)f_n\int_0^Tp(s)(T-s)^{\alpha-2}E_{\alpha,\alpha-1}(-\lambda_n(T-s)^{\alpha})ds\bigg|\le c_{48}
\|p\|_{L^{\infty}(0,T)}\|f\|_{D(A^{\frac{1}{\alpha}})}.
\end{equation}
Hence we see that $\partial_t^{\alpha}u\in L^{\infty}(0,T;H)$, the series \eqref{Eq78.0} is convergent in $L^{\infty}(0,T;H)$ and 
\begin{equation*}
    \|\partial_t^{\alpha}u\|_{L^{\infty}(0,T;H)}\le c_{49}\|p\|_{L^{\infty}(0,T)}.
\end{equation*}
So, we get the first inequality in \eqref{EqTh4.1}.

Since the series is convergent in $L^{\infty}(0,T;H)$, we have
\begin{equation*}
    \partial_t^{\alpha}h(t)=\Phi[f]p(t)-\sum\limits_{n=1}^{\infty}\lambda_nf_n \Big(p(t)\star t^{\alpha-1}E_{\alpha,\alpha}(-\lambda_nt^{\alpha})\Big)\Phi[e_n]\qquad\qquad\qquad\qquad\qquad
\end{equation*}
\begin{equation*}
    +\sum\limits_{n=1}^{\infty}\lambda_n\mathcal{P}_n^{(2)}(t)\cdot f_n \Big(p(r)\star r^{\alpha-1}E_{\alpha,\alpha}(-\lambda_nr^{\alpha})\Big)\Big|_{r=T}\Phi[e_n]
\end{equation*}
\begin{equation}\label{Eq79}
    \quad+\sum\limits_{n=1}^{\infty}\lambda_n\mathcal{P}_n^{(3)}(t)\cdot f_n \Big(p(r)\star r^{\alpha-2}E_{\alpha,\alpha-1}(-\lambda_nr^{\alpha})\Big)\Big|_{r=T}\Phi[e_n],
\end{equation}
for $0<t<T$. Setting
\begin{equation*}
    \mathcal{K}_1(t)=-\sum\limits_{n=1}^{\infty}\lambda_nf_n  t^{\alpha-1}E_{\alpha,\alpha}(-\lambda_nt^{\alpha})\Phi[e_n],\quad \mathcal{K}_2(t)=\sum\limits_{n=1}^{\infty}\lambda_n\mathcal{P}_n^{(2)}(t)\cdot f_n  t^{\alpha-1}E_{\alpha,\alpha}(-\lambda_nt^{\alpha})\Phi[e_n]
\end{equation*}
and by (C1) and \eqref{Eq25}, we can see that $\mathcal{K}_1\in C[0,T]$ and $\mathcal{K}_2\in C[0,T]$.  Moreover, we set
\begin{equation*}
    \mathcal{K}_3(t)=\sum\limits_{n=1}^{\infty}\lambda_n\mathcal{P}_n^{(3)}(t)\cdot f_n E_{\alpha,\alpha}(-\lambda_nt^{\alpha})\Phi[e_n].
\end{equation*}
The same reason as above and by \eqref{Eq74}, we have $\mathcal{K}_3\in C[0,T]$. Therefore
\begin{equation*}
    \partial_t^{\alpha}h(t)=p(t)\Phi[f]+\int_0^t\mathcal{K}_1(t-s)p(s)ds\qquad\qquad\qquad\qquad\qquad\qquad\qquad
\end{equation*}
\begin{equation*}
    \qquad\qquad+\int_0^T\mathcal{K}_2(T-s)p(s)ds+\int_0^T(T-s)^{\alpha-2}\mathcal{K}_3(T-s)p(s)ds,\quad 0<t<T,
\end{equation*}
that is,
\begin{equation*}
    p(t)=\frac{\partial_t^{\alpha}h(t)}{\Phi[f]}-\frac{1}{\Phi[f]}\int_0^t\mathcal{K}_1(t-s)p(s)ds\qquad\qquad\qquad\qquad\qquad\qquad\qquad\qquad\qquad\qquad\qquad
\end{equation*}
\begin{equation*}
    -\frac{1}{\Phi[f]}\int_0^T\mathcal{K}_2(T-s)p(s)ds -\frac{1}{\Phi[f]}\int_0^T(T-s)^{\alpha-2}\mathcal{K}_3(T-s)p(s)ds,\quad 0<t<T
\end{equation*}
by (C2). Hence
\begin{equation*}
    |p(t)|\le c_{50}\|\partial_t^{\alpha}h\|_{L^{\infty}(0,T)}+c_{51}\|\mathcal{K}_1\|_{C[0,T]}\int_0^t|p(s)|ds\qquad\qquad\qquad\qquad\qquad\qquad\qquad\qquad\quad
\end{equation*}
\begin{equation*}
    +c_{52}\|\mathcal{K}_2\|_{C[0,T]}\int_0^T|p(s)|ds+c_{53}\|\mathcal{K}_2\|_{C[0,T]}\int_0^T(T-s)^{\alpha-2}|p(s)|ds,\quad 0<t<T.
\end{equation*}
Applying an inequality of Vollter-Fredholm type weakly singular kernel $(T-s)^{\alpha-2}$ (e.g., Theorem 4, \cite{Zeng}), we see
\begin{equation*}
    |p(t)|\le c_{54}\|\partial_t^{\alpha}h\|_{L^{\infty}(0,T)},\quad 0<t<T,
\end{equation*}
that is, the second inequality in \eqref{EqTh4.1} is proved. Thus, the proof of Theorem 4.1 is complete.

\section{CONCLUSION}

In this paper, we have investigated a time-fractional diffusion-wave equation involving a positive self-adjoint operator with compact inverse in a separable Hilbert space. We first established the well-posedness of the direct problem with non-local initial conditions, providing a spectral representation of the solution. Two inverse problems were then studied in detail. In Inverse Problem 1, we considered the identification of a time-independent source term using additional information at an intermediate time. We proved the uniqueness of the source and presented an explicit reconstruction method based on eigenfunction expansion. In Inverse Problem 2, we addressed the recovery of a time-dependent coefficient in the source term from scalar output data. By reducing the problem to a Volterra integral equation of the second kind, we developed a solvability framework for the unknown coefficient. These results contribute to the theory of inverse problems for fractional differential equations and may have further applications in anomalous diffusion modeling and memory-dependent processes. 

An open problem remains concerning the behaviour of solutions in certain cases of the parameter $\gamma.$ As noted in Remark 2.1, in such cases, there is no available information about an upper bound for the largest zero 
$\eta_N$, which poses analytical challenges in establishing uniform estimates or further spectral bounds. This issue points to a potential direction for future research.

\section{ACKNOWLEDGEMENTS}
The author would like to thank the referees for their careful reading and valuable suggestions on this paper.

\section{CONFLICTS OF INTEREST}
This work does not have any conflicts of interest.

\end{document}